\documentclass[12pt]{amsart}
\usepackage[a4paper,margin=27mm]{geometry}
\usepackage[]{fontenc}
\usepackage{lmodern,microtype}
\usepackage{amsmath,amssymb,amsthm,mathtools}
\usepackage[numbers,sort&compress]{natbib}
\usepackage[colorlinks=true,linkcolor=blue,citecolor=blue,urlcolor=blue,pagebackref=true]{hyperref}
\numberwithin{equation}{section}
\newtheorem{theorem}{Theorem}[section]
\newtheorem{lemma}[theorem]{Lemma}
\newtheorem{proposition}[theorem]{Proposition}

\theoremstyle{definition}
\newtheorem{definition}[theorem]{Definition}
\theoremstyle{remark}
\newtheorem{remark}[theorem]{Remark}
\newcommand{\tr}{\operatorname{tr}}
\newcommand{\A}{\mathcal A}
\newcommand{\B}{\mathcal B}
\newcommand{\C}{\mathcal C}
\newcommand{\D}{\mathcal D}
\newcommand{\F}{\mathfrak F}
\newcommand{\ii}{\sqrt{-1}}

\allowdisplaybreaks[1]
\title{The Space of $k$-Positive Potentials}
\begin{document}

\author{Xishen Jin}
	\address{Xishen Jin\\ Department of Mathematics\\ Renmin University of China \\ Beijing\\ 100872\\ China\\} \email{jinxishen@ruc.edu.cn}
    \subjclass[2020]{Primary 32Q15; Secondary 35B65, 35J70, 58E10.}

\keywords{Weak geodesics, $k$-positive potentials, complex Hessian equations}
\maketitle
\begin{abstract}
In this paper, we study weak geodesics in the space of $k$-positive potentials on a
compact K\"ahler manifold. We prove that any two smooth strictly
$k$-positive potentials can be joined by a weak geodesic of class $C^{1,1}$. As an application, we prove that the pair $(\mathcal H_k,d_2)$ is a metric space.
\end{abstract}

\section{Introduction}
Let $(X,\omega)$ be a compact K\"ahler manifold of complex dimension
$n$, without boundary. Mabuchi~\cite{Mabuchi} defined a Riemannian metric on the space of K\"ahler metrics, under which it becomes a non-positive curved infinite dimensional symmetric space. Later, Semmes~\cite{Semmes} and
Donaldson~\cite{Donaldson} related its geodesic equation to a homogeneous
complex Monge--Amp\`ere equation on a manifold (with boundary) of one dimension higher. This relation gives a way to study
geodesics by elliptic equations on a product manifold. Chen~\cite{Chen}
proved the existence of weak geodesics and obtained a bound for the
complex Hessian. The bound for all real second derivatives was proved
by Chu--Tosatti--Weinkove~\cite{CTW}.

Complex Hessian equations give a natural extension of this setting. For a real vector $\lambda=(\lambda_1,\ldots,
\lambda_n)$, let $\sigma_j(\lambda)$ denote its $j$-th elementary
symmetric polynomial. For $1\leq k\leq n$, define
\[
 \Gamma_k=\{\lambda\in\mathbb R^n:
                 \sigma_j(\lambda)>0\text{ for }1\leq j\leq k\}.
\]
We use the same notation for the cone of Hermitian matrices with
eigenvalues in $\Gamma_k$. Given $\varphi\in C^\infty(X,\mathbb R)$,
write $\omega_\varphi=\omega+\ii\partial\bar\partial\varphi$.
The space of smooth $k$-positive potentials is
\[
 \mathcal H_k=\{\varphi\in C^\infty(X,\mathbb R):
          \lambda_\omega(\omega_\varphi)\in\Gamma_k\}.
\]
Here $\lambda_\omega$ denotes the eigenvalues with respect to $\omega$.

Following Mabuchi \cite{Mabuchi}, Guo--Smith--Tong~\cite{GST} studied the formal Riemannian metric
\begin{equation}\label{eq:metric}
 \langle v_1,v_2\rangle_\varphi
       =\int_X v_1v_2\,\omega_\varphi^k\wedge\omega^{n-k}
\end{equation}
on $\mathcal H_k$. They derived its connection and geodesic equation.
Let $u=u(x,t)$ be a smooth path in $\mathcal H_k$.
In local holomorphic coordinates $(z^1,\ldots,z^n)$, write
$\omega=\ii g_{p\bar q}\,dz^p\wedge d\bar z^q$ and set
\[
 G_u^{p\bar q}
 =k\frac{\ii\,dz^p\wedge d\bar z^q\wedge\omega_u^{k-1}
                 \wedge\omega^{n-k}}
              {\omega_u^k\wedge\omega^{n-k}}.
\]
The coefficients \(G_u^{p\bar q}\) arise from the linearization $L_u v=\sum\limits_{p,q=1}^n G_u^{p\bar q}v_{p\bar q}$ of
\[
u\longmapsto \log\frac{\omega_u^k\wedge\omega^{n-k}}{\omega^n}.
\]
Following~\cite[Section 2.3]{GST}, the geodesic equation is
\[
 u_{tt}-\sum_{p,q=1}^n G_u^{p\bar q}u_{tp}u_{t\bar q}=0,
 \qquad
 u_{tp}=\frac{\partial^2u}{\partial t\,\partial z^p}.
\]
Multiplying by $\omega_u^k\wedge\omega^{n-k}$ gives
\begin{equation}\label{eq:geodesic}
 u_{tt}\,\omega_u^k\wedge\omega^{n-k}
 -k \ii\,\partial_X u_t\wedge\bar\partial_X u_t
       \wedge\omega_u^{k-1}\wedge\omega^{n-k}=0
\end{equation}
where $\partial_X$ and $\bar\partial_X$ denote spatial
derivatives. We
introduce the complex variable $\tau=t+\ii s$ and extend $u$ to $u(\cdot,\tau)=u(\cdot,t)$ 
independently of $s$. Let
\[
 \Sigma=\{\tau\in\mathbb C:0<\operatorname{Re}\tau<1\},
 \qquad
 \pi:X\times\Sigma\longrightarrow X
\]
be the projection. On $X\times\Sigma$, set
\[
 \Omega_u=\pi^*\omega+\ii\partial\bar\partial u,
\]
where $\partial$ and $\bar\partial$ now act on both the
spatial and time variables. Since $u$ is independent of $s$,
\[
 u_{\tau\bar\tau}=\frac14u_{tt},
 \qquad
 u_{p\bar\tau}=\frac12u_{pt}.
\]
Expanding the wedge product gives
\[
 \begin{split}
&\Omega_u^{k+1}\wedge(\pi^*\omega)^{n-k}\\
 =
 &\frac{k+1}{4}\,\ii\,d\tau\wedge d\bar\tau\wedge
 \left(
 u_{tt}\,\omega_u^k\wedge\omega^{n-k}
 -
 k \ii\,\partial_Xu_t\wedge\bar\partial_Xu_t
 \wedge\omega_u^{k-1}\wedge\omega^{n-k}
 \right).
 \end{split}
\]
Thus the geodesic equation with endpoints
$\varphi_0,\varphi_1\in\mathcal H_k$ becomes the degenerate
Dirichlet problem
\begin{equation}\label{eq:geodesic-dirichlet}
 \begin{cases}
 \bigl(\pi^*\omega+\ii\partial\bar\partial u\bigr)^{k+1}
       \wedge(\pi^*\omega)^{n-k}=0
       & \text{on }X\times\Sigma,\\
 u(x,\tau)=\varphi_0(x)
       & \text{when }\operatorname{Re}\tau=0,\\
 u(x,\tau)=\varphi_1(x)
       & \text{when }\operatorname{Re}\tau=1.
 \end{cases}
\end{equation}
This equation is degenerate: its right-hand side vanishes,
and the background form $\pi^*\omega$ has no component
in the time direction.
The corresponding complex operator was also introduced by
Chen--He~\cite[Section 3]{ChenHeClass}.
For $k=n$, equation~\eqref{eq:geodesic} is the usual K\"ahler
geodesic equation. For $k=1$, a fixed rescaling gives Donaldson's
equation for geodesics in the space of volume forms~\cite{DonaldsonVolume}.
Chen--He~\cite{ChenHe} constructed weak geodesics in that space.
Chu~\cite{ChuVolume} proved their full $C^{1,1}$ regularity.
For $1<k<n$, the admissible cone is larger than the positive cone,
and its ellipticity can degenerate in additional directions.

The spatial estimates for complex Hessian equations are central to
this problem. Hou--Ma--Wu~\cite{HMW} obtained a second order estimate
in terms of the gradient. Dinew--Ko\l odziej~\cite{DK} proved a
Liouville theorem which gives the gradient estimate by rescaling.
Chu--McCleerey~\cite{CM} established real Hessian estimates for a broad
class of degenerate equations in complex geometry. On the potential
theory side, \AA hag--Czy\.{z}--Lu--Rashkovskii~\cite{ACLR} studied
geodesic envelopes and a minimum principle for complex Hessian
equations with a semipositive background form.

In this paper, we seek solutions to the Dirichlet problem \eqref{eq:geodesic-dirichlet} in the admissible class defined in Section~\ref{sec:cone}. A main result of the paper is the following $C^{1,1}$ estimate in both real space and time variables. 

\begin{theorem}\label{thm:main}
Let $\varphi_0,\varphi_1\in\mathcal H_k$. The weak geodesic $U$ joining
them belongs to $C^{1,1}(X\times[0,1])$. More precisely,
\[
 \|U\|_{C^1(X\times[0,1])}
       +\|\nabla^2_{X,t}U\|_{L^\infty(X\times[0,1])}\leq C,
\]
where the norms are respect to the fixed product metric, and $C$ depends only
on $(X,\omega)$ and the endpoint data $\varphi_0$, $\varphi_1$. In particular, the slices $\omega_{\varphi_t}$ are weakly $k$-positive.
Moreover, $U$ is also
the unique continuous viscosity solution of \eqref{eq:geodesic-dirichlet} with these boundary values.
\end{theorem}

We next consider the path length induced by \eqref{eq:metric}.
For a piecewise smooth path $\gamma:[0,1]\to\mathcal H_k$, set
\[
 L_2(\gamma)=\int_0^1
 \left(\int_X|\dot\gamma_s|^2\,
       \omega_{\gamma_s}^k\wedge\omega^{n-k}\right)^{1/2}ds.
\]
Let $d_2(\varphi_0,\varphi_1)$ be the infimum of $L_2(\gamma)$
over all such paths joining $\varphi_0$ to $\varphi_1$.

\begin{theorem}\label{thm:hessian-distance}
The pair $(\mathcal H_k,d_2)$ is a metric space.
If $U$ is the weak geodesic joining $\varphi_0$ and $\varphi_1$
in Theorem~\ref{thm:main}, then for every $t\in[0,1]$,
\begin{equation}\label{eq:hessian-distance}
 d_2(\varphi_0,\varphi_1)^2
 =\int_X |U_t(\cdot,t)|^2\,
          \omega_{U(\cdot,t)}^k\wedge\omega^{n-k}.
\end{equation}
\end{theorem}

The paper is organized as follows.
In Section~\ref{sec:cone}, we describe the admissible cone,
introduce the regularized equation, and obtain uniform bounds for $u$ and $u_t$.
In Section~\ref{sec:spatial-estimates}, by taking the partial Legendre
transform in $t$, we estimate the gradient and the complex Hessian
of the transformed potential $Z$ along $X$.
In Section~\ref{sec:real-hessian}, we estimate the mixed derivatives
on the boundary, bound $u_{tt}$, and control the real Hessian of $Z$
along $X$. In Section~\ref{mx:sec:mixed}, we bound the mixed derivative
quotient that appears in the inverse Legendre formula and obtain a
uniform bound for all real second derivatives of $u$.
We use these estimates in Section~\ref{sec:existence} to solve
the regularized Dirichlet problem by the continuity method.
In Section~\ref{sec:limit}, we prove comparison for the limiting
equation and Theorem \ref{thm:main}.
Finally, in Section~\ref{sec:hessian-distance}, we prove Theorem \ref{thm:hessian-distance}.

\medskip
\noindent\textbf{Note on this version.}
The results of this paper on the existence of weak geodesics and
the induced path distance overlap in part with those in the first
version of \cite{GWX}, where Ghosh, Wang, and Xu studied a more
general class of complex Hessian equations. We have since combined
the two works into a joint paper. The present preprint will not be
submitted separately for publication.

\section{The Admissible Cone And Preliminary Estimates}
\label{sec:cone}

In this section, we specify the admissible branch of the
geodesic equation and introduce a perturbed Dirichlet problem.
We then obtain bounds for the solution and its first time
derivative. In particular, these bounds are independent of the perturbation
parameter.

\subsection{The Admissible Cone}

Let $u$ be a smooth real function on $X\times[0,1]$.
At a fixed point, choose a unitary frame with
respect to $\omega$. Write
\[
 A_u=I+(u_{i\bar j}),\qquad
 b_u=(u_{it}),\qquad
 c_u=u_{tt},
\]
where $I$ is the identity matrix and spatial derivatives
are covariant derivatives for the fixed K\"ahler metric.
We regard $b_u$ as a column vector and set
\[
 B_u=
 \begin{pmatrix}
 A_u&b_u\\
 b_u^*&c_u
 \end{pmatrix}.
\]
Here the superscript $*$ denotes the conjugate transpose.

For a Hermitian matrix $A$ and a Hermitian variation $H$,
write
\[
 D\sigma_k(A)[H]
 =
 \left.\frac{d}{dr}\right|_{r=0}\sigma_k(A+rH).
\]
This derivative can also be written as
\[
 D\sigma_k(A)[H]
 =
 \operatorname{tr}\bigl(T_{k-1}(A)H\bigr),
\]
where $T_{k-1}(A)$ is the Newton matrix defined by
\[
 T_{k-1}(A)
 =
 \sum_{j=0}^{k-1}(-1)^j
 \sigma_{k-1-j}(A)A^j,
 \qquad \sigma_0=1.
\]
In a unitary frame diagonalizing $A$, this matrix is
diagonal, with entries
\[
 \bigl(T_{k-1}(A)\bigr)_{i\bar i}
 =
 \sigma_{k-1}(\lambda|i).
\]
Here $\lambda$ is the vector of eigenvalues of $A$,
and $\lambda|i$ means deleting
$\lambda_i$ from $\lambda=(\lambda_1,\ldots,\lambda_n)$.
And $D\sigma_k(A)[H]$ is given by
\[
 D\sigma_k(A)[H]
 =
 \sum_{i=1}^n
 \sigma_{k-1}(\lambda|i)\,H_{i\bar i}.
\]
In particular, for a column vector
$b=(b_1,\ldots,b_n)^T$ expressed in the same frame,
\[
 D\sigma_k(A)[bb^*]
 =
 \sum_{i=1}^n
 \sigma_{k-1}(\lambda|i)\,|b_i|^2,
\]
where $b^*$ denotes the conjugate transpose of $b$.
For a block matrix
\[
 B=\begin{pmatrix}A&b\\b^*&c\end{pmatrix},
\]
define
\[
 P(B)=c\,\sigma_k(A)-D\sigma_k(A)[bb^*].
\]
Then the geodesic equation~\eqref{eq:geodesic} is equivalent to $P(B_u)=0$.
Indeed, in the chosen frame,
\[
 \omega_u^k\wedge\omega^{n-k}
 =
 \binom{n}{k}^{-1}\sigma_k(A_u)\,\omega^n.
\]
Differentiating this identity in the direction of the
Hermitian matrix $b_ub_u^*$ gives
\[
 k \ii\,\partial_Xu_t\wedge\bar\partial_Xu_t
       \wedge\omega_u^{k-1}\wedge\omega^{n-k}
 =
 \binom{n}{k}^{-1}
 D\sigma_k(A_u)[b_ub_u^*]\,\omega^n.
\]
Substitution into~\eqref{eq:geodesic} proves the assertion.
This block operator is the complex operator considered
in~\cite[Section 3]{ChenHeClass}.

We first record an elementary identity.

\begin{lemma}\label{lem:rankone}
For every Hermitian matrix $A$, every $b\in\mathbb C^n$,
and every $r\in\mathbb R$,
\[
 \sigma_k(A+rbb^*)
 =
 \sigma_k(A)+r\,D\sigma_k(A)[bb^*].
\]
In particular, if $c>0$, then
\begin{equation}\label{eq:schurpoly}
 P(B)=c\,\sigma_k\left(A-\frac{bb^*}{c}\right).
\end{equation}
\end{lemma}

\begin{proof}
The real symmetric version is proved in
\cite[Lemma 2.5]{GS}.
The same argument applies to Hermitian matrices.
Indeed, if $b\ne0$, we choose a unitary frame in which
$b=|b|e_1$. Then $A+rbb^*$ differs from $A$ only
in its first diagonal entry. Hence $\sigma_k(A+rbb^*)$ is affine in $r$.
Its derivative at $r=0$ is $D\sigma_k(A)[bb^*]$.

Taking $r=-1/c$ and multiplying by $c$, we get
the second identity.
\end{proof}

We first recall the following standard properties of $\Gamma_k$;
see~\cite{Garding,CNS}.
\begin{proposition}\label{prop:garding-cone}
The cone $\Gamma_k$ has the following properties:
\begin{enumerate}
 \item $\Gamma_k$ is open and convex. It is the connected
 component of $\{A:\sigma_k(A)>0\}$ containing the
 positive definite Hermitian matrices.

 \item If $A\in\Gamma_k$ and $H\geq0$, then
 $A+H\in\Gamma_k$.

 \item If $A\in\Gamma_k$, then $T_{k-1}(A)$ is positive
 definite. Consequently,
 \[
  D\sigma_k(A)[H]
  =\operatorname{tr}\bigl(T_{k-1}(A)H\bigr)\geq0
  \qquad\text{for }H\geq0,
 \]
 and the inequality is strict when $H\ne0$.
\end{enumerate}
\end{proposition}

We define
\begin{equation}\label{eq:cone}
 \mathcal C_k
 =
 \left\{
 \begin{pmatrix}A&b\\b^*&c\end{pmatrix}:
 c>0,\quad A-\frac{bb^*}{c}\in\Gamma_k
 \right\}.
\end{equation}
Then we have the following equivalent characterization for $\mathcal C_k$.

\begin{lemma}\label{lem:cone-equivalence}
For a Hermitian block matrix $B=\begin{pmatrix}A&b\\b^*&c\end{pmatrix}$, the following
conditions are equivalent:
\[
 A\in\Gamma_k,\qquad P(B)>0;
\]
and
\[
 c>0,\qquad A-\frac{bb^*}{c}\in\Gamma_k.
\]
\end{lemma}

\begin{proof}
Assume first that $A\in\Gamma_k$ and $P(B)>0$.
Since $bb^*\geq0$, the third item of Proposition \ref{prop:garding-cone} implies
\[
 c\,\sigma_k(A)
 >
 D\sigma_k(A)[bb^*]\geq0.
\]
As $\sigma_k(A)>0$, this implies $c>0$.

For $0\leq r\leq1$, set
\[
 A_r=A-\frac{rbb^*}{c}.
\]
Lemma~\ref{lem:rankone} gives
\[
 \begin{aligned}
 \sigma_k(A_r)
 &=
 \sigma_k(A)-\frac r cD\sigma_k(A)[bb^*]\\
 &=
 (1-r)\sigma_k(A)+\frac r cP(B)>0.
 \end{aligned}
\]
Thus $A_r$ is a continuous path in $\{\sigma_k>0\}$.
It starts at $A\in\Gamma_k$, so it stays in the same
connected component according to Proposition \ref{prop:garding-cone}. In particular,
\[
 A_1=A-\frac{bb^*}{c}\in\Gamma_k.
\]

Conversely, suppose that $c>0$ and
$M=A-bb^*/c\in\Gamma_k$.
Since $bb^*/c\geq0$, applying the second item of Proposition \ref{prop:garding-cone}, we have
\[
 A=M+\frac{bb^*}{c}\in\Gamma_k.
\]
Equation~\eqref{eq:schurpoly} then gives
$P(B)=c\,\sigma_k(M)>0$.
\end{proof}

We call a smooth function $u$ strictly admissible if
$B_u=
 \begin{pmatrix}
 A_u&b_u\\
 b_u^*&c_u
 \end{pmatrix}\in\mathcal C_k$ at every point.
This condition is independent of the chosen spatial
unitary frame. In particular, a strictly admissible
function satisfies $u_{tt}>0$ and $A_u\in\Gamma_k$.

\begin{lemma}\label{lem:cone-convexity}
The set $\mathcal C_k$ is an open convex cone.
If $B\in\mathcal C_k$ and $H\geq0$, then
$B+H\in\mathcal C_k$.
\end{lemma}

\begin{remark}
The real symmetric version of the convexity statement
was proved in~\cite[Theorem 2.1]{HXZ} using hyperbolic
polynomials.
Here we give a short and absolutely different proof in the Hermitian setting.
\end{remark}

\begin{proof}
Openness follows directly from the definition.
For every $a>0$,
\[
 aA-\frac{(ab)(ab)^*}{ac}
 =
 a\left(A-\frac{bb^*}{c}\right).
\]
Since $\Gamma_k$ is a cone, $aB\in\mathcal C_k$
whenever $B\in\mathcal C_k$.

We next prove convexity.
Let $B_0,B_1\in\mathcal C_k$, and write their blocks
as $A_j,b_j,c_j$, where $j=0,1$.
For $0\leq\theta\leq1$, let
\[
 B_\theta=(1-\theta)B_0+\theta B_1.
\]
For each $\xi\in\mathbb C^n$, the Cauchy--Schwarz
inequality gives
\[
 \frac{|(1-\theta)b_0^*\xi+\theta b_1^*\xi|^2}
      {(1-\theta)c_0+\theta c_1}
 \leq
 (1-\theta)\frac{|b_0^*\xi|^2}{c_0}
 +\theta\frac{|b_1^*\xi|^2}{c_1}.
\]
Hence, in the order of Hermitian matrices,
\[
 \frac{b_\theta b_\theta^*}{c_\theta}
 \leq
 (1-\theta)\frac{b_0b_0^*}{c_0}
 +\theta\frac{b_1b_1^*}{c_1}.
\]
Writing $M_j=A_j-b_jb_j^*/c_j\in \Gamma_k$, we obtain
\[
 M_\theta=A_\theta-\frac{b_\theta b_\theta^*}{c_\theta}
 \geq (1-\theta)M_0+\theta M_1.
\]
The right-hand side belongs to $\Gamma_k$ by convexity.
According to Proposition \ref{prop:garding-cone}, $M_\theta$ also belongs to $\Gamma_k$ and $P(B_\theta)>0$.

$A_\theta\in \Gamma_k$ since $\Gamma_k$ is convex. Thus $B_\theta\in\mathcal C_k$.

Finally, for $c>0$, completion of the square
\[
\begin{aligned}
\begin{pmatrix}\xi\\z\end{pmatrix}^{\!*}
B
\begin{pmatrix}\xi\\z\end{pmatrix}
&=\xi^*A\xi+2\operatorname{Re}(\xi^*bz)+c|z|^2\\
&=\xi^*\left(A-\frac{bb^*}{c}\right)\xi
+c\left|z+\frac{b^*\xi}{c}\right|^2
\end{aligned}
\]
gives
\[
 \xi^*\left(A-\frac{bb^*}{c}\right)\xi
 =
 \min_{z\in\mathbb C}
 \begin{pmatrix}\xi\\z\end{pmatrix}^{\!*}
 B
 \begin{pmatrix}\xi\\z\end{pmatrix}.
\]
We note that the matrix $A-bb^*/c$ is the Schur complement
of the entry $c$ in $B$. If $H\geq0$, adding $H$ increases the expression being
minimized. The bottom-right entry of $B+H$ is still positive.
Its Schur complement is therefore no smaller than that
of $B$. The monotonicity of $\Gamma_k$ now gives
$B+H\in\mathcal C_k$.
\end{proof}

\subsection{The Perturbed Equation}

Set
\[
 f(M)=\sigma_k(M)^{1/k},
 \qquad M\in\Gamma_k.
\]
The function $f$ is positive, smooth, homogeneous of
degree one, and concave on $\Gamma_k$; see~\cite{Garding,CNS}. For a Hermitian variation $H$, its derivative is
\[
 Df(M)[H]
 =
 \frac{1}{k}\sigma_k(M)^{\frac1k-1}
 \operatorname{tr}\bigl(T_{k-1}(M)H\bigr).
\]
Thus under the trace
pairing, $Df(M)$ is represented by
\[
 \frac{1}{k}\sigma_k(M)^{\frac1k-1}T_{k-1}(M).
\]
This matrix is positive definite because
$\sigma_k(M)>0$ and $T_{k-1}(M)>0$ on $\Gamma_k$.

For $B\in\mathcal C_k$, we define
\begin{equation}\label{eq:rootoperator}
 \mathfrak F(B)
 =
 \left[
 c\,f\left(A-\frac{bb^*}{c}\right)
 \right]^{1/2}.
\end{equation}
The square root in this definition is useful for
concavity.

\begin{lemma}\label{lem:rootoperator}
The function $\mathfrak F$ is smooth, concave, and
homogeneous of degree one on $\mathcal C_k$.
Its derivative satisfies
\[
 D\mathfrak F(B)[H]>0
 \qquad\text{for every } H\geq0,\quad H\ne0.
\]
\end{lemma}

\begin{proof}
Smoothness follows from $c>0$ and $f>0$.
Since $f(aM)=af(M)$ for $a>0$, we have
$\mathfrak F(aB)=a\mathfrak F(B)$.

To prove concavity, we use the notation from the proof of
Lemma~\ref{lem:cone-convexity}. The proof there gives
\[
 M_\theta
 \geq (1-\theta)M_0+\theta M_1.
\]
The monotonicity and concavity of $f$ imply
\[
 f(M_\theta)
 \geq
 f((1-\theta)M_0+\theta M_1)
 \geq
 (1-\theta)f(M_0)+\theta f(M_1).
\]
Consequently, we have
\[
 \begin{aligned}
 \mathfrak F(B_\theta)
 &\geq
 \sqrt{
 ((1-\theta)c_0+\theta c_1)
 ((1-\theta)f(M_0)+\theta f(M_1))
 }\\
 &\geq
 (1-\theta)\sqrt{c_0f(M_0)}
 +\theta\sqrt{c_1f(M_1)}
 \end{aligned}
\]
where we use Cauchy--Schwarz in the last inequality.
This proves concavity.

We next compute the derivative along $H$.
Write
\[
 H=\begin{pmatrix}E&d\\d^*&h\end{pmatrix},
 \qquad
 M=A-\frac{bb^*}{c}.
\]
Differentiation gives
\[
 \delta M
 =
 E-\frac{db^*+bd^*}{c}+\frac{hbb^*}{c^2}
 =
 \begin{pmatrix}I&-b/c\end{pmatrix}
 H
 \begin{pmatrix}I\\-b^*/c\end{pmatrix}.
\]
Thus
\begin{equation}\label{eq:root-linearization}
 2\mathfrak F(B)D\mathfrak F(B)[H]
 =
 h f(M)+c\,Df(M)[\delta M].
\end{equation}
If $H\geq0$, then $h\geq0$ and $\delta M\geq0$.
Both terms on the right-hand side are nonnegative.

Suppose in addition that $H\ne0$.
If $h>0$, the first term is positive.
If $h=0$, we consider the principal submatrix 
\[
 \begin{pmatrix}
 E_{i\bar i}&d_i\\
 \overline{d_i}&h
 \end{pmatrix}
\]
for every $i=1,2,\cdots,n$. Since $H\geq0$, this submatrix is nonnegative.
Its determinant therefore satisfies
\[
 E_{i\bar i}h-|d_i|^2\geq0.
\]
As $h=0$, we obtain $d_i=0$ for every $i$.
Thus $\delta M=E\geq0$ and $E\ne0$.
Since $Df(M)$ is positive definite,
$Df(M)[E]>0$. This proves the strict inequality.
\end{proof}

In particular, $u\mapsto\mathfrak F(B_u)$ is elliptic
on strictly admissible functions. For $0<\varepsilon\leq1$, we consider the following perturbed equation
\begin{equation}\label{eq:root}
 \begin{cases}
 \displaystyle
 u_{tt}\,
 f\left(
 I+(u_{i\bar j})
 -\dfrac{(u_{it}u_{t\bar j})}{u_{tt}}
 \right)=\varepsilon
 &\text{on }X\times(0,1),\\[6pt]
 u(\cdot,0)=\varphi_0,\qquad
 u(\cdot,1)=\varphi_1,
 \end{cases}
\end{equation}
with $B_u\in\mathcal C_k$.
Equivalently,
\[
 \mathfrak F(B_u)=\sqrt{\varepsilon}.
\]
Using~\eqref{eq:schurpoly}, the equation also gives
\[
 P(B_u)=\varepsilon^k u_{tt}^{\,1-k}.
\]

From now on, we assume that $u$ is a smooth strictly admissible
solution of~\eqref{eq:root}.

\subsection{$C^0$ And First Time Derivative Estimates}\label{sbsec:C0}

For the $C^0$ estimate, we construct a lower barrier with the prescribed
boundary values. Set
\[
 h(x,t)=(1-t)\varphi_0(x)+t\varphi_1(x),
 \qquad
 d(x)=\varphi_1(x)-\varphi_0(x).
\]
By the definition of $h$, we have
\[
 A_h=(1-t)A_{\varphi_0}+tA_{\varphi_1}.
\]
Since $\varphi_0$ and $\varphi_1$ are strictly $k$-positive,
the convexity of $\Gamma_k$ gives $A_h\in\Gamma_k$.

The matrices $A_h$ form a compact family in $\Gamma_k$.
Indeed, $X\times[0,1]$ is compact and the eigenvalues
depend continuously on $(x,t)$.
We can therefore choose $\delta>0$ such that
\[
 A_h-\delta I\in\Gamma_k
 \qquad\text{on }X\times[0,1].
\]
Since $f$ is continuous and positive on this compact
family, the number
\[
 m=\min_{X\times[0,1]}f(A_h-\delta I)
\]
is positive.

For a constant $K>0$ to be chosen, we define
\[
 \underline u=h-Kt(1-t).
\]
Its associated Hermitian block matrix is 
\[
B_{\underline u}=\begin{pmatrix}
 A_{\underline u}&b_{\underline u}\\
 b_{\underline u}^*&c_{\underline u}
\end{pmatrix}\]
where
\[
 A_{\underline u}=\omega+\ii\partial_X\bar\partial_X {\underline u}=A_h,\qquad
 b_{\underline u}=\partial_X {\underline u}_t=\partial_Xd,\qquad
 c_{\underline u}={\underline u}_{tt}=2K.
\]
Choose $K$ large enough that
\[
 \frac{|\partial_Xd|_\omega^2}{2K}\leq\delta
 \quad\text{on }X,
 \qquad
 2Km>1.
\]
Then the Schur complement
\[
 \begin{aligned}
 M_{\underline u}
 &:=
 A_{\underline u}
 -\frac{b_{\underline u}b_{\underline u}^*}
        {c_{\underline u}}\\
 &=
 A_h-\frac{(\partial_Xd)(\partial_Xd)^*}{2K}
 \geq A_h-\delta I.
 \end{aligned}
\]
It follows that $M_{\underline u}\in\Gamma_k$ and
$f(M_{\underline u})\geq m$.
Thus $\underline u$ is strictly admissible and
\[
 c_{\underline u}f(M_{\underline u})
 \geq2Km>1\geq\varepsilon.
\]
Also,
\[
 \mathfrak F(B_{\underline u})=\sqrt{c_{\underline u}f(M_{\underline u})}>\sqrt{\varepsilon} ,\qquad\underline u(\cdot,0)=\varphi_0,\qquad
 \underline u(\cdot,1)=\varphi_1.
\]

\begin{proposition}\label{prop:basic-estimates}
There is a constant $C$, depending only on $\varphi_0$, $\varphi_1$ and $(X,\omega)$, such that
every smooth strictly admissible solution of
\eqref{eq:root} satisfies
\begin{equation}\label{eq:low}
 \|u\|_{L^\infty(X\times[0,1])}
 +
 \|u_t\|_{L^\infty(X\times[0,1])}
 \leq C.
\end{equation}
In particular, the constant is independent of $0<\varepsilon\leq1$.
\end{proposition}

\begin{proof}
We first prove $\underline u\leq u$.
Suppose that $\underline u-u$ has a positive maximum.
Since $u=\underline u$ at $t=0,1$, this maximum
occurs at an interior point.

At that point its real Hessian is nonpositive.
The associated Hermitian block matrix is also
nonpositive, i.e.
\[
 B_{\underline u}\leq B_u.
\]

By Lemma~\ref{lem:cone-convexity}, the segment $B_{tu+(1-t)\underline u}$ from $B_{\underline u}$
 to $B_u$ stays in $\mathcal C_k$.
Integrating along this segment and applying Lemma~\ref{lem:rootoperator},we have
\[\mathfrak F(B_u)
-\mathfrak F(B_{\underline u})=\int_0^1 D\mathfrak F(B_{tu+(1-t)\underline u})[B_u-B_{\underline u}]dt\geq 0.\]
Hence, $
 \mathfrak F(B_{\underline u})
 \leq\mathfrak F(B_u)=\sqrt{\varepsilon}$. 
This contradicts the strict barrier inequality.
Therefore $\underline u\leq u$.

Since $u_{tt}>0$, the function $t\mapsto u(x,t)$ is
convex for each fixed $x$.
We have
\begin{equation}\label{eq:barrier-bounds}
 u\leq h
\end{equation}
where $h$ is the line joint $u(x,0)$ and $u(x,1)$, by the definition. In particular, if
\[
 C_0=\max\{\|\varphi_0\|_{L^\infty(X)},
           \|\varphi_1\|_{L^\infty(X)}\},
\]
then
\[
 -C_0-\frac K4\leq u\leq C_0.
\]

We next estimate $u_t$ at $t=0,1$.
At $t=0$, both
$u-\underline u$ and $h-u$ vanish and are nonnegative
for $t\geq0$.
Since
\[
 h_t=d,\qquad
 \underline u_t=d-K+2Kt,
\]
we obtain
\[
 d-K\leq u_t(x,0)\leq d.
\]
Similarly, at $t=1$
\[
 d\leq u_t(x,1)\leq d+K.
\]
Finally, $u_{tt}>0$ implies that $u_t(x,t)$ is
increasing in $t$. Therefore
\[
 d(x)-K
 \leq u_t(x,0)
 \leq u_t(x,t)
 \leq u_t(x,1)
 \leq d(x)+K.
\]
Consequently,
\[
 \|u_t\|_{L^\infty(X\times[0,1])}
 \leq \|\varphi_1-\varphi_0\|_{L^\infty(X)}+K.
\]
Combining this with the bound for $u$ proves
\eqref{eq:low}.
\end{proof}

In the following sections, $C$ denotes a positive
constant that may change from line to line.
Unless stated otherwise, it depends only on the fixed
background, the functions $\varphi_0$ and $\varphi_1$. In particular, it is independent
of $\varepsilon$.

\section{Gradient And Complex Hessian Estimate In The $X$ Direction}
\label{sec:spatial-estimates}

In this section, we use a partial Legendre transform
in the time variable to obtain gradient and complex
Hessian estimates in the $X$ directions.

\subsection{The Partial Legendre Transform}
\label{sp:sec:legendre}

Partial Legendre transforms have been used to study
nonlinear elliptic equations; see \cite{GP,DR}. We first apply the transform in the time variable
to rewrite~\eqref{eq:root}.

Let $u$ be a smooth strictly admissible solution of
\eqref{eq:root}.
We define the partial Legendre transform of $u$ in the time variable by
\begin{equation}\label{sp:eq:extension}
 \widehat Z(x,s)
 =
 \inf_{0\leq t\leq1}\{u(x,t)-st\},
 \qquad (x,s)\in X\times\mathbb R.
\end{equation}
Since $u_{tt}>0$, the minimizer is unique.
Set
\[
 a_0(x)=u_t(x,0),\qquad a_1(x)=u_t(x,1),
\]
and define
\[
 \Omega=\{(x,s)\in X\times\mathbb R:
                  a_0(x)<s<a_1(x)\}.
\]
For $(x,s)\in\Omega$, the minimizer of \eqref{sp:eq:extension} lies in $(0,1)$
and is determined by
\[
 s=u_t(x,t).
\]
The implicit function theorem gives a smooth function
$t=t(x,s)$ on $\Omega$.
We write $Z$ for the restriction of $\widehat Z$
to $\Omega$. Thus
\[
 Z(x,s)=u(x,t(x,s))-s\,t(x,s).
\]

For the standard differentiation formulas for partial
Legendre transforms, see~\cite[Section 2.2]{GP}.
Differentiating
$s=u_t(x,t)$ and $Z=u-st$, we have
\[
 t_s=\frac1{u_{tt}},\qquad
 t_i=-\frac{u_{it}}{u_{tt}},\qquad
 Z_s=-t,\qquad Z_i=u_i.
\]
Differentiating once more, we obtain
\begin{equation}\label{sp:eq:Legendre-identities}
 Z_{ss}=-\frac1{u_{tt}},\qquad
 Z_{is}=\frac{u_{it}}{u_{tt}},\qquad
 Z_{i\bar j}
 =
 u_{i\bar j}
 -\frac{u_{it}u_{t\bar j}}{u_{tt}}.
\end{equation}

For convenience, we write
\[
 A=A_u,\qquad b=b_u,\qquad c=c_u=u_{tt},
\]
We also set
\[
 W=-Z_{ss}=\frac1c,\qquad
 p=\partial_XZ_s=\frac bc,\qquad
 M=I+(Z_{i\bar j})=A-\frac{bb^*}{c}.
\]
Strict admissibility gives $M\in\Gamma_k$ and $W>0$.

With notation above, Equation~\eqref{eq:root} now becomes
\begin{equation}\label{sp:eq:ZPDE}
 f(M)+\varepsilon Z_{ss}=0.
\end{equation}
By~\eqref{eq:low}, the functions $a_0,a_1$ are uniformly
bounded. Thus $\Omega$ lies in a fixed bounded interval
of the $s$ variable. Since $Z=u-st$, we also have
\[
 \|Z\|_{L^\infty(\Omega)}\leq C.
\]

We next record the boundary bounds needed for the
complex Hessian estimate.

\begin{lemma}\label{sp:lem:boundarySchur}
At the boundary $t=j$, where $j=0,1$,
let $A_j=A_{\varphi_j}$. Then at $t=j$, we have
\[
 M=A_j-\frac{bb^*}{c},
 \qquad
 \frac{|b|^2}{c}
 =\operatorname{tr}A_j-\operatorname{tr}M
 \leq\operatorname{tr}A_j.
\]
Consequently, the complex Hessian and the real
Hessian $\nabla_X^2 Z$ of $Z$ in the space direction are uniformly bounded on $\partial\Omega$.
Moreover,
\[
 \frac{|p|^2}{W}\leq C
 \qquad\text{on }\partial\Omega.
\]
\end{lemma}

\begin{proof}
The first identity follows from
\eqref{sp:eq:Legendre-identities} and $u(x,j)=\varphi_j (x)$ at $t=j$.
Taking trace, we obtain
\[
 \frac{|b|^2}{c}
 =\operatorname{tr}A_j-\operatorname{tr}M\leq \operatorname{tr}A_j\leq C
\]
at the boundary, since $M\in\Gamma_k\subset\Gamma_1$.

We assume that the eigenvalues of $M$
are $\mu_1\geq\cdots\geq\mu_n$. Then $M\leq A_j$ implies $\mu_1\leq C$. Since $M\in\Gamma_k\subset\Gamma_1$, we also have
\[\mu_n>-\sum_{i=1}^{n-1}\mu_i
 \geq-(n-1)\mu_1\geq C.\]
Hence the complex Hessian of $Z$ is bounded.

For the real Hessian, differentiating the Legendre
transform gives
\[
 \nabla_X^2Z
 =
 \nabla_X^2\varphi_j
 -
 \frac{d_Xu_t\otimes d_Xu_t}{c} \leq \nabla_X^2\varphi_j.
\]
The last term is bounded since $|b|^2/c\leq C$ where $b=\partial_X u_t$ at the boundary.
Since $\frac{|p|^2}{W}=\frac{|b|^2}{c}$, the last inequality is obvious.
\end{proof}

\subsection{The Complex Hessian Estimate}
\label{sp:sec:complex}

To obtain the complex Hessian estimate, we first establish a uniform estimate of Hou--Ma--Wu type~\cite{HMW} by the maximum principle.

\begin{proposition}\label{sp:prop:HMW}
There are constants $\varepsilon_0>0$ and $C$,
depending only on $(X,\omega)$ and $\varphi_0,\varphi_1$, such that
\begin{equation}\label{sp:eq:HMW}
\sup_\Omega|\partial_X\bar\partial_X Z|
 \leq C\left(1+\sup_\Omega|\partial_XZ|^2\right)
\end{equation}
for $0<\varepsilon\leq\varepsilon_0$.
\end{proposition}

\begin{proof}
Write
\[
 P=|\partial_X Z|^2,\qquad K=1+\sup_\Omega P.
\]
Choose a fixed constant $L>1+\sup_\Omega|Z|$ and set
\[
 \phi(P)=-\frac12\log\left(1-\frac{P}{2K}\right),
 \qquad
 \psi(Z)=-A_0\log\left(1+\frac{Z}{2L}\right),
 \qquad
 A_0=3LD,
\]
where $D\geq1$ will be chosen below. For $(x,s)\in\overline\Omega$
and a unit vector $\xi\in T_x^{1,0}X$, consider the auxiliary function as in \cite{HMW},
\[
 W(x,s,\xi)
 =M(\xi,\bar\xi)\exp\bigl(\phi(P)+\psi(Z)\bigr).
\]
Direct differentiation gives
\begin{equation}\label{sp:eq:weights}
 \frac1{4K}\leq\phi'\leq\frac1{2K},
 \qquad
 \phi''=2(\phi')^2,
 \qquad
 D\leq-\psi'\leq3D,
 \qquad
 \psi''=\frac{(\psi')^2}{A_0}.
\end{equation}
For fixed $L$ and $D$, both $\phi$ and $\psi$ have bounded
oscillation, independently of $K$ and $\varepsilon$.

Let $(x_0,s_0,\xi_0)$ be a maximum point of $W$. Since the
exponential factor is independent of $\xi$, the vector $\xi_0$
is an eigenvector of $M(x_0,s_0)$ corresponding to its largest
eigenvalue. If $(x_0,s_0)\in\partial\Omega$, the estimate follows
from Lemma~\ref{sp:lem:boundarySchur} and the bounded oscillations
of $\phi$ and $\psi$. We may therefore assume that
$(x_0,s_0)\in\Omega$.

Choose holomorphic normal coordinates at $x_0$ in which
$M(x_0,s_0)$ is diagonal, with eigenvalues
$\lambda_1\geq\cdots\geq\lambda_n$ and
$\xi_0=\partial/\partial z^1$ at $x_0$. Set
\[
 h=\frac{M_{1\bar1}}{g_{1\bar1}}.
\]
Then $h\leq\lambda_1$ nearby, with equality at $(x_0,s_0)$.
Since $\operatorname{tr}M>0$, we have $\lambda_1>0$ at this point.
Thus $h$ is positive in a neighborhood of $(x_0,s_0)$, and the
smooth function
\[
 Q=\log h+\phi(P)+\psi(Z)
\]
has a local maximum there. All calculations below are made at
$(x_0,s_0)$ unless otherwise indicated.

The linearized operator of~\eqref{sp:eq:ZPDE} is
\[
 \mathcal D
 =f^{i\bar j}\nabla_i\nabla_{\bar j}
  +\varepsilon\partial_{ss},
\]
where $f^{i\bar j}$ denotes the derivative of $f$ with respect
to $M_{i\bar j}$. Set
\[
 F_i=\frac{\partial f}{\partial\lambda_i},
 \qquad
 \mathcal A=\sum_iF_i,
 \qquad
 \mathcal B=-D^2f.
\]
The form 
\begin{equation}\label{eq:B-spec}
\mathcal B(H,H^*)
=
-\sum_{p,q}\frac{\partial^2 f}{\partial \lambda_p \partial \lambda_q}H_{p\bar p}\overline{H_{q\bar q}}
+\sum_{p\ne q}
\frac{F_q-F_p}{\lambda_p-\lambda_q}|H_{p\bar q}|^2,
\end{equation} is nonnegative by the concavity of $f$.
When applied to derivatives in complex directions, $\mathcal B$
denotes its complex bilinear extension.

The following inequalities will be used:
\begin{equation}\label{sp:eq:elementary}
 \lambda_1>0,\quad
 |\lambda_i|\leq(n-1)\lambda_1\quad(i>1),
 \quad
 F_1\leq\cdots\leq F_n,
 \quad
 \mathcal A\geq f(I)=\binom nk^{1/k}\geq1.
\end{equation}
The first two inequalities follow from $\operatorname{tr}M>0$.
The third follows from the symmetry and concavity of $f$.
For the last inequality, concavity and homogeneity give
\[
 \mathcal A=Df(M)[I]
 =f(M)+Df(M)[I-M]
 \geq f(I)=\binom nk^{1/k}\geq1.
\]

Put
\[
 a_i=\frac{h_i}{h},
 \qquad
 a_s=\frac{h_s}{h}.
\]
Since $Q$ has a local maximum at $(x_0,s_0)$, its first
derivatives vanish. Hence
\begin{equation}\label{sp:eq:critical}
 a_i=-\phi'P_i-\psi'Z_i,
 \qquad
 a_s=-\phi'P_s-\psi'Z_s.
\end{equation}
Moreover, $\mathcal DQ\leq0$. To compute $\mathcal DQ$, we first
compute $\mathcal DZ$, $\mathcal DP$, and $\mathcal Dh$.

Since $f$ is homogeneous of degree one and $M=I+(Z_{i\bar j})$
in the chosen coordinates, the equation gives
\begin{equation}\label{sp:eq:DZ}
 \mathcal DZ
 =Df(M)[M-I]+\varepsilon Z_{ss}
 =f(M)-\mathcal A+\varepsilon Z_{ss}
 =-\mathcal A.
\end{equation}

Differentiating~\eqref{sp:eq:ZPDE} in the $k$ direction, we obtain
\[
 \sum_i F_i Z_{i\bar i k}+\varepsilon Z_{ssk}=0.
\]
The K\"ahler commutation formulas give
\[
 Z_{k\bar i i}=Z_{i\bar i k},
 \qquad
 Z_{ki\bar i}
 =Z_{i\bar i k}+\sum_p R_{k\bar i i}{}^p Z_p.
\]
In particular,
\[
 \sum_i F_i Z_{ki\bar i}
 -\sum_{i,p}F_i R_{k\bar i i}{}^p Z_p
 +\varepsilon Z_{ssk}=0.
\]
Applying the product rule to $P=|\partial_X Z|^2$ and using these
identities, we find
\begin{align}
 \mathcal DP
 ={}&\sum_{i,p}F_i\bigl(|Z_{pi}|^2+|Z_{p\bar i}|^2\bigr)
     +2\varepsilon|\partial_X Z_s|^2 \notag\\
 &+2\operatorname{Re}\sum_k Z_{\bar k}
       \left(\sum_iF_i Z_{i\bar i k}+\varepsilon Z_{ssk}\right)
   +\sum_{i,k,p}F_i R_{k\bar i i}{}^p Z_pZ_{\bar k} \notag\\
 \geq{}&\sum_{i,p}F_i\bigl(|Z_{pi}|^2+|Z_{p\bar i}|^2\bigr)
     +2\varepsilon|\partial_X Z_s|^2-C_0P\mathcal A.
 \label{sp:eq:DP}
\end{align}
Here the term in parentheses vanishes by the differentiated
equation, and $C_0$ depends only on the background metric.

Next, differentiate~\eqref{sp:eq:ZPDE} in the $1,\bar1$ directions.
Commuting the spatial derivatives as in \cite[(2.20)]{HMW}, we obtain
\begin{equation}\label{sp:eq:Dlambda}
 \begin{aligned}
 \mathcal Dh
 &=\mathcal B(M_1,M_{\bar1})
   +\sum_i F_i(\lambda_1-\lambda_i)R_{1\bar1 i\bar i}\\
 &\geq\mathcal B(M_1,M_{\bar1})-C_0\lambda_1\mathcal A.
 \end{aligned}
\end{equation}
Here $M_1$ and $M_{\bar1}$ denote covariant derivatives of $M$ and the constant $C_0$ depends only on the background curvature bound.

Since $Q=\log h+\phi(P)+\psi(Z)$ has a local maximum at
$(x_0,s_0)$, the maximum principle gives
\begin{align*}
 0\geq\mathcal DQ
 ={}&\frac{\mathcal Dh}{h}
     +\phi'\mathcal DP+\psi'\mathcal DZ\\
 &+\sum_i F_i\left(
      -|a_i|^2+\phi''|P_i|^2+\psi''|Z_i|^2
   \right)\\
 &+\varepsilon\left(
      -a_s^2+\phi''P_s^2+\psi''Z_s^2
   \right)
\end{align*}
where $a_i$ and $a_s$ are defined in \eqref{sp:eq:critical}. 

We first estimate the last line.
By the $s$-derivative condition in~\eqref{sp:eq:critical},
\[
 a_s^2=(\phi'P_s+\psi'Z_s)^2
 \leq 2(\phi')^2P_s^2+2(\psi')^2Z_s^2.
\]
Since $\phi''=2(\phi')^2$ and $\phi',\psi''\geq0$, we obtain
\begin{equation}\label{sp:eq:timeerror}
 \begin{aligned}
 &\varepsilon\left(
   -a_s^2+\phi''P_s^2+\psi''Z_s^2
 \right)\\
 &\qquad\geq -2\varepsilon(\psi')^2Z_s^2
 \geq -18\varepsilon D^2,
 \end{aligned}
\end{equation}
where the last inequality follows from $|\psi'|\leq3D$ and
$|Z_s|\leq1$.

Since $Z_{i\bar i}=\lambda_i-1$ at $(x_0,s_0)$,
\eqref{sp:eq:weights} and $(\lambda_i-1)^2\geq\lambda_i^2/2-1$
give
\[
 \phi'\sum_{i,p}F_i|Z_{p\bar i}|^2
 \geq\frac1{8K}\sum_iF_i\lambda_i^2-\phi'\mathcal A.
\]
Combining this inequality with~\eqref{sp:eq:DZ},
\eqref{sp:eq:DP}, \eqref{sp:eq:Dlambda},
and~\eqref{sp:eq:timeerror}, we obtain
\begin{align}
 0\geq{}&
 \frac{\mathcal B(M_1,M_{\bar1})}{\lambda_1}
 -\sum_iF_i|a_i|^2
 +\phi''\sum_iF_i|P_i|^2
 +\psi''\sum_iF_i|Z_i|^2
 \notag\\
 &+\frac1{8K}\sum_iF_i\lambda_i^2
 +(-\psi'-C_*)\mathcal A
 -18\varepsilon D^2
 \label{sp:eq:master}
\end{align}
where the nonnegative terms involving $|Z_{pi}|^2$ and $|\partial_X Z_s|^2$ is dropped.
Since $P\leq K$ and $\phi'\leq1/(2K)$, the constant $C_*$
is independent of $D$, $K$, and $\varepsilon$. Choose
\[
 D=\max\{C_*+2,1\},
 \qquad
 \theta=\frac1{1+2A_0}.
\]
Then $0<\theta<1/2$ and $-\psi'-C_*\geq2$.

\medskip
\noindent
\emph{Case 1: $\lambda_n<-\theta\lambda_1$.}

By~\eqref{sp:eq:critical} and $\phi''=2(\phi')^2$, we have
\[
 -\sum_iF_i|a_i|^2+\phi''\sum_iF_i|P_i|^2
 \geq-2(\psi')^2\sum_iF_i|Z_i|^2
 \geq-18D^2K\mathcal A.
\]
Since $F_n\geq\mathcal A/n$, inequality~\eqref{sp:eq:master}
implies
\[
 0\geq
 \frac{\theta^2\lambda_1^2}{8nK}\mathcal A
 -18D^2K\mathcal A-18\varepsilon D^2.
\]
For $\varepsilon\leq1$, using $K,\mathcal A\geq1$, we obtain
\[
 \lambda_1^2\leq 8nK\theta^{-2}\mathcal A^{-1}(18D^2K\mathcal A+18D^2)\leq288nD^2\theta^{-2}K^2.
\]

\medskip
\noindent
\emph{Case 2: $\lambda_n\geq-\theta\lambda_1$.}

Let
\[
 I=\{i:F_i>\theta^{-1}F_1\}.
\]
Then $1\notin I$. For the indices outside $I$,
\eqref{sp:eq:critical} gives
\begin{equation}\label{sp:eq:smallweights}
 -\sum_{i\notin I}F_i|a_i|^2
 +\phi''\sum_{i\notin I}F_i|P_i|^2
 \geq -2(\psi')^2\sum_{i\notin I}F_i|Z_i|^2\geq-18D^2\theta^{-1}KF_1.
\end{equation}

For $i\in I$, we have $F_i>F_1$ and hence $\lambda_i<\lambda_1$.
At the chosen point, $h_i=Z_{1\bar1 i}=Z_{i\bar1 1}$. Since $f$ is concave, we have
\begin{align}
 \frac{\mathcal B(M_1,M_{\bar1})}{\lambda_1}
 &\geq
 \sum_{i\in I}
 \frac{F_i-F_1}{\lambda_1(\lambda_1-\lambda_i)}
 |Z_{1\bar1 i}|^2
 \notag\\
 &=\sum_{i\in I} \frac{1-F_1/F_i}{1-\lambda_i/\lambda_1}F_i|a_i|^2\geq \frac{1-\theta}{1+\theta}\sum_{i\in I}F_i|a_i|^2\\
 &\geq
 (1-2\theta)\sum_{i\in I}F_i|a_i|^2.
 \label{sp:eq:strongweights}
\end{align}

The first derivative equations and $\phi''=2(\phi')^2$ also give
\[
 \phi''|P_i|^2
 =2|a_i+\psi'Z_i|^2
 \geq2\theta|a_i|^2
 -\frac{2\theta}{1-\theta}(\psi')^2|Z_i|^2=2\theta|a_i|^2
 -\psi''|Z_i|^2.
\]
Since
\[
 \frac{2\theta}{1-\theta}=\frac1{A_0},
 \qquad
 \psi''=\frac{(\psi')^2}{A_0},
\]
the last term is cancelled by $\psi''|Z_i|^2$.
Combining this estimate with~\eqref{sp:eq:strongweights}, we obtain
\[
 \frac{\mathcal B(M_1,M_{\bar1})}{\lambda_1}
 +\sum_{i\in I}F_i
 \left(-|a_i|^2+\phi''|P_i|^2+\psi''|Z_i|^2\right)
 \geq 0.
\]

If $\lambda_1^2\geq 288D^2\theta^{-1}K^2$, then
\[
 \frac{1}{16K}\sum_i F_i\lambda_i^2
 \geq \frac{F_1\lambda_1^2}{16K}
 \geq 18D^2\theta^{-1}KF_1.
\]
Thus half of the quadratic term in~\eqref{sp:eq:master} absorbs
the negative term in~\eqref{sp:eq:smallweights}, leaving
\[
 0\geq
 \frac{1}{16K}\sum_i F_i\lambda_i^2
 +2\mathcal A-18\varepsilon D^2.
\]
This is impossible for $\varepsilon\leq(18D^2)^{-1}$, since
$\mathcal A\geq1$. Thus $\lambda_1\leq CK$ in this case as well.

We have proved $\lambda_1(x_0,s_0)\leq CK$. By the maximality of
$W$ and the bounded oscillations of $\phi$ and $\psi$, the same
bound, with a possibly larger constant, holds throughout $\Omega$.
The remaining eigenvalues are controlled by~\eqref{sp:eq:elementary}.
This proves~\eqref{sp:eq:HMW}.
\end{proof}

\subsection{The Gradient Estimate}
\label{sp:sec:gradient}

For the blow-up argument below, we use $\widehat Z$ from
\eqref{sp:eq:extension}, since $\widehat Z(\cdot,s)$ is defined
on all of $X$ for each fixed $s$. Since $u$ is strictly convex in $t$,
\[
 \widehat Z(x,s)=
 \begin{cases}
  \varphi_0(x), & s\leq a_0(x),\\
  Z(x,s),       & a_0(x)<s<a_1(x),\\
  \varphi_1(x)-s, & s\geq a_1(x).
 \end{cases}
\]
Since $Z_i=u_i$ and $Z_s=-t$ on $\Omega$, we have
\[
 (Z_i,Z_s)=
 \begin{cases}
  ((\varphi_0)_i,0),  & s=a_0(x),\\
  ((\varphi_1)_i,-1), & s=a_1(x),
 \end{cases}
\]
where the derivatives are taken from within $\Omega$.
These are precisely the first derivatives of the two exterior
branches, so $\widehat Z$ is $C^1$ across both interfaces.

For each fixed $\varepsilon$, $u_{tt}$ has
a positive minimum (may depending on $\varepsilon$) on $X\times[0,1]$. The Legendre derivative
formulas therefore give bounded second derivatives of $Z$ on
$\Omega$. The two exterior branches are smooth, and their first
derivatives agree with those of $Z$ along the interfaces. It follows
that $\widehat Z\in C^{1,1}_{\mathrm{loc}}(X\times\mathbb R)$,
with a bound that may depend on $\varepsilon$.

For each fixed $x$, the function $\widehat Z(x,\cdot)$ is concave
because it is the infimum of affine functions of $s$. Moreover,
$\widehat Z_s=-t$ on $\Omega$ and takes the values $0$ and $-1$
on the two exterior branches. Thus
\[
 -1\leq\widehat Z_s\leq0.
\]
Taking $t=0$ and $t=1$ in the defining infimum of $\widehat Z$ also gives
\begin{equation}\label{sp:eq:obstacles}
 \widehat Z\leq\varphi_0,
 \qquad
 \widehat Z\leq\varphi_1-s.
\end{equation}

For each fixed $s$, by locality of weak derivatives and the
admissibility of $Z$ on $\Omega$ and of the endpoint data,
\[
 I+(\widehat Z_{i\bar j})\in\overline{\Gamma_k}
 \qquad\text{almost everywhere on }X.
\]
Thus $\widehat Z(\cdot,s)$ is $k$-subharmonic relative to $\omega$.
The complex Hessian bound in Proposition~\ref{sp:prop:HMW}
also holds almost everywhere for $\widehat Z(\cdot,s)$, after
enlarging the constant to account for the endpoint bounds.

\begin{theorem}\label{sp:thm:gradient}
For all sufficiently small $\varepsilon$, every smooth
strictly admissible solution of~\eqref{eq:root} satisfies
\[
 \sup_\Omega|\partial_XZ|
 +
 \sup_\Omega|\partial_X\bar\partial_XZ|
 \leq C,
 \qquad
 \sup_{X\times[0,1]}|\partial_Xu|\leq C.
\]
The same bounds hold for $\widehat Z$
almost everywhere.
\end{theorem}

\begin{proof}
Suppose that the gradient bound fails.
Then there is a sequence of solutions $u_j$, with parameters
$\varepsilon_j$, such that its partial Legendre transform $\widehat Z_j$ satisfies
\[
 G_j=\sup_{X\times\mathbb R}
             |\partial_X\widehat Z_j|
 \longrightarrow\infty.
\]
Suppose that $(x_j,s_j)$ are the maximum points.
Outside the domain \[
\Omega_j=\{(x,s):u_{j,t}(x,0)<s<u_{j,t}(x,1)\}
\] 
of $Z_j$, the gradient $\partial_X\widehat Z_j$
equals either $\partial_X\varphi_0$ or $\partial_X\varphi_1$.
Since $G_j\to\infty$, the maximum points $(x_j,s_j)$ lie in
$\Omega_j$ for all large $j$. The uniform bound for $u_{j,t}$ then implies that $s_j$ lies
in a fixed compact interval.

After passing to a subsequence, we may assume that
$(x_j,s_j)\to(x_\infty,s_\infty)$. Starting from fixed holomorphic
coordinates near $x_\infty$, we translate the origin to $x_j$
and make a complex linear change of coordinates so that the metric
at the origin is the identity. Since these linear changes have
uniformly bounded inverses, the resulting coordinate domains
contain a common ball $B_\rho(0)\subset\mathbb C^n$ for all large $j$.

Writing $\widehat Z_j$ in these coordinates, we define
\begin{equation}\label{sp:eq:scaling}
 V_j(z,r)=\widehat Z_j\left(
 \frac{z}{G_j},\,s_j+\frac{\sqrt{\varepsilon_j}}{G_j}r
 \right),
 \qquad (z,r)\in B_{\rho G_j}(0)\times\mathbb R.
\end{equation}
Since $G_j\to\infty$, these domains eventually contain every
compact subset of $\mathbb C^n\times\mathbb R$.

On each fixed compact set, we have
\[
 |\partial_zV_j|\leq C,
 \qquad
 |(V_j)_r|\leq\frac{\sqrt{\varepsilon_j}}{G_j}
 \longrightarrow0.
\]
By Proposition~\ref{sp:prop:HMW} and its extension to
$\widehat Z_j$, we also have
\[
 |\partial_z\bar\partial_z V_j|
 \leq \frac{C(1+G_j^2)}{G_j^2}
 \leq C'
 \qquad\text{almost everywhere}.
\]
Indeed, each derivative with respect to $z$ or $\bar z$
contributes a factor of $G_j^{-1}$.

On the slice $r=0$, the local bounds for $V_j$ and $\Delta_z V_j$
give uniform interior $W^{2,p}$ bounds for every finite $p$.
Taking $p>2n$ and $0<\alpha<1-2n/p$, we obtain, after passing
to a subsequence,
\[
 V_j(\cdot,0)\longrightarrow w
 \qquad\text{in }C^{1,\alpha}_{\mathrm{loc}}(\mathbb C^n).
\]
Moreover,
\[
 |V_j(z,r)-V_j(z,0)|
 \leq |r|\frac{\sqrt{\varepsilon_j}}{G_j},
\]
so $V_j$ converges locally uniformly on $\mathbb C^n\times\mathbb R$
to the same limit $w(z)$, which is independent of $r$.
By the normalization $|\partial_z V_j(0,0)|=1$ and convergence
in $C^1$ on the slice $r=0$, we have
\[
 |\partial w(0)|=1.
\]

Since the sequence $s_j$ lies in a fixed compact interval,
\eqref{eq:low} and the defining formula for $\widehat Z_j$
give a uniform bound for $\widehat Z_j(\cdot,s_j)$ on $X$.
It follows that $w$ is bounded on $\mathbb C^n$. Moreover,
the local $C^{1,\alpha}$
convergence on the slice $r=0$ and the uniform bound of $|\partial_X V_j|$ imply that $w$ is globally Lipschitz.

On the slice $r=0$, the cone condition for $\widehat Z_j$ gives
\[
 \frac{1}{G_j^2}\omega(\frac{z}{G_j})
 +\sqrt{-1}\,\partial_z\bar\partial_z V_j(\cdot,0)
 \in\overline{\Gamma_k}
\]
almost everywhere.
On compact sets, $G_j^{-2}\omega$ tends uniformly to zero. Since the
complex Hessians of $V_j(\cdot,0)$ are locally uniformly bounded, $w$ is $k$-subharmonic with respect to
the Euclidean metric after passing to the limit.

In the variables $(z,r)$ of~\eqref{sp:eq:scaling}, the functions
$\varphi_0(x)$ and $\varphi_1(x)-s$ converge locally uniformly to
the constants $d_0=\varphi_0(x_\infty)$ and $d_1=\varphi_1(x_\infty)-s_\infty$,
respectively. Then
\eqref{sp:eq:obstacles} gives $w\leq\min\{d_0,d_1\}$. Since $w$ is subharmonic
and $|\partial w(0)|=1$, the strong maximum principle gives
\[
 w<\min\{d_0,d_1\}\qquad\text{on }\mathbb C^n.
\]

On each compact subset $D\Subset \mathbb C^n\times\mathbb R$, there exists $\eta>0$ such that
\[w\leq \min\{d_0,d_1\}-\eta.\]
By local uniform convergence, for all sufficiently
large $j$
\[V_j\leq \min\{d_0,d_1\}-\frac\eta2<\min\{d_0,d_1\}\]
on $D$. Since $\widehat Z_j$ equals $\varphi_0$ or $\varphi_1-s$
outside $\Omega_j$, the compact set $D$ lies in the rescaled domain of
$Z_j$. Thus the rescaled equation holds there for all sufficiently
large $j$, i.e.
\begin{equation}\label{sp:eq:scaledPDE}
 f\left(
 G_j^{-2}I+
 g_j^{-1}\partial_z\bar\partial_zV_j
 \right)
 +(V_j)_{rr}=0,
 \qquad g_j\longrightarrow I.
\end{equation}
Here $g_j$ is the background metric evaluated in the
rescaled coordinates. 
The factor $\sqrt{\varepsilon_j}$ in
\eqref{sp:eq:scaling} makes the coefficient of
$(V_j)_{rr}$ equal to one.

We claim that no smooth function with strictly
$k$-positive complex Hessian can touch $w$ from below.
Suppose that a function $\chi(z)$ does so at $z_0$.
By subtracting a sufficiently small multiple of
$|z-z_0|^2$, we may make the contact strict while
keeping its complex Hessian strictly $k$-positive
near $z_0$.
Choose $\delta>0$ small enough that
\[
 f(\partial\bar\partial\chi)>2\delta
\]
on a small neighborhood.
Then $\chi(z)-\delta r^2$ touches $w(z)$ strictly
from below at $(z_0,0)$.

Set $\eta(z,r)=\chi(z)-\delta r^2$. By strict contact,
$w-\eta$ has a unique minimum at $(z_0,0)$ on a small
closed ball $\overline B$ centered there. Since $V_j\to w$
uniformly on $\overline B$, the minimum of $V_j-\eta$ on
$\overline B$ is attained at an interior point
$(z_j,r_j)\to(z_0,0)$ for all sufficiently large $j$.
Setting $c_j=(V_j-\eta)(z_j,r_j)$, we see that $\eta+c_j$
touches $V_j$ from below at $(z_j,r_j)$.
At these contact points, ellipticity and
\eqref{sp:eq:scaledPDE} imply
\[
 f\left(
 G_j^{-2}I+
 g_j^{-1}\partial\bar\partial\chi
 \right)-2\delta\leq0.
\]
Passing to the limit contradicts the choice of $\delta$.
This proves the claim.

Next, we prove that $w$ is maximal $k$-subharmonic function on $\mathbb C^n$, i.e. for every $k$-subharmonic function $v$ and any compact subset $D\Subset \mathbb C^n$, the inequality $w\geq v$ on $\mathbb C^n\setminus D$ implies $w\geq v$ on $D$.
We argue by contradiction. Suppose that $v>w$ at some inner point of $D$.
Choose $R$ such that $\overline D\subset\{|z|<R\}$.
For sufficiently small $\delta>0$, the function
\[
 v_\delta=v+\delta(|z|^2-R^2)
\]
satisfies $v_\delta>w$ at some point in $D$. Since
$v_\delta<v$ on $\partial D$, we still have
$v_\delta<w$ on $\partial D$.

By replacing $v$ with $\max\{v,-N\}$ for sufficiently
large $N$, we may assume that $v$ is bounded since as a $k$-subharmonic function, $v$ is bounded above.
Choose a subdomain $D'\Subset D$ containing a point $z_*$ with
$v_\delta(z_*)>w(z_*)$, and such that $v_\delta<w$ near $\partial D'$.
Let $\rho_\eta$ be a standard radial mollifier and set
\[
 v_{\delta,\eta}
 =v*\rho_\eta+\delta(|z|^2-R^2)
 \qquad\text{on }\overline{D'}.
\]
Convolution preserves $k$-subharmonicity, so
\[
 \partial\bar\partial v_{\delta,\eta}
 =\partial\bar\partial(v*\rho_\eta)+\delta I\in\Gamma_k.
\]
Since $v$ is subharmonic, $v*\rho_\eta\geq v$. Thus
$v_{\delta,\eta}(z_*)>w(z_*)$. For sufficiently small $\eta$,
the strict inequality near $\partial D'$ also gives
$v_{\delta,\eta}<w$ on $\partial D'$.

Hence
\[
 m_\eta=\max_{\overline{D'}}(v_{\delta,\eta}-w)>0
\]
is attained in $D'$. At a maximum point, $v_{\delta,\eta}-m_\eta$
touches $w$ from below and has strictly $k$-positive complex
Hessian. This contradicts the preceding claim. Therefore
$v\leq w$ on $D$.

Thus $w$ is maximal in the sense of complex Hessian
potential theory.
The Liouville theorem of
Dinew--Ko{\l}odziej~\cite[Theorem 3.2]{DK} applies:
a bounded maximal $k$-subharmonic function on
$\mathbb C^n$ with bounded gradient is constant.
This contradicts $|\partial w(0)|=1$.

The gradient $\partial_X \widehat Z$ is therefore uniformly bounded.
Proposition~\ref{sp:prop:HMW} then gives the uniform
complex Hessian bound.
Finally, the identity $Z_i=u_i$ implies
\[
 |\partial_Xu(x,t)|
 =
 |\partial_XZ(x,u_t(x,t))|.
\]
This proves the gradient bound for $u$ and completes
the proof.
\end{proof}

\section{The Uniform Real Hessian Estimate}\label{sec:real-hessian}

In this section, we first bound the mixed derivatives of $u$ on
the boundary and $u_{tt}$ on $X\times[0,1]$. We then prove a
uniform bound for the real Hessian $\nabla_X^2 Z$ on $\Omega$.

\subsection{Boundary Mixed Derivatives And The Estimate For $u_{tt}$}
\label{sp:sec:time}

\begin{proposition}\label{sp:prop:boundarymixed}
For all sufficiently small $\varepsilon>0$, every smooth strictly
admissible solution of~\eqref{eq:root} satisfies
\[
 |\nabla_Xu_t|_{t=0,1}\leq C,
 \qquad
 0<u_{tt}\leq C\quad\text{on }X\times[0,1],
\]
where $C$ is independent of $\varepsilon$. In particular,
$W=-Z_{ss}\geq C^{-1}$ on $\Omega$.
\end{proposition}

\begin{proof}
Let $\mathcal L$ be the linearization of $u\mapsto\F(B_u)$ at $u$,
where $\F$ is the concave operator in Lemma~\ref{lem:rootoperator}.
We write
\[
 \mathcal Lv=a^{\alpha\beta}\nabla_\alpha\nabla_\beta v,
\]
where $\alpha,\beta$ range over the real coordinate indices on
$X\times[0,1]$ and $(a^{\alpha\beta})$ is positive definite.

Let $\underline u$ be the strict subsolution constructed in Subsection \ref{sbsec:C0}
and set $w=u-\underline u\geq0$.
Choose $\eta>0$ such that $B_{\underline u}-\eta I\in\C_k$ everywhere.
For sufficiently small $\varepsilon$, we have
$\F(B_{\underline u}-\eta I)\geq\F(B_u)=\sqrt{\varepsilon}$.
By concavity, there exists a fixed $\eta'>0$
\begin{equation}\label{sp:eq:strictbarrier}
 \begin{aligned}
 \mathcal L(w)
 &=D\F(B_u)[B_{\underline u}-B_u]\\
 &\geq \F(B_{\underline u}-\eta I)-\F(B_u)
       +\eta\,\tr D\F(B_u)\\
 &\geq \eta'\tr a
 \end{aligned}
\end{equation}
where $a=(a^{\alpha\beta})$. Here we used
$\tr a\leq C\tr D\F(B_u)$, since $a$ is obtained from
$D\F(B_u)$ by a fixed linear map.

Set $h=(1-t)\varphi_0+t\varphi_1$ and $v=d_Xu-d_Xh$.
Then $v=0$ at $t=0,1$. Differentiating the equation along $X$ gives
\[
 a^{\alpha\beta}\nabla_i(\nabla^2u)_{\alpha\beta}=0.
\]
Since $du$ is uniformly
bounded by~\eqref{eq:low} and Theorem~\ref{sp:thm:gradient},
and $h$ is fixed, commuting the third derivatives in $\mathcal L(d_X u)$, we obtain
\begin{equation}\label{sp:eq:covector}
 |\mathcal Lv|\leq C\tr a.
\end{equation}

For $\tau>0$, we set $\rho_\tau=(|v|^2+\tau^2)^{1/2}-\tau$.
Then we have
\[
 \mathcal L\rho_\tau
 \geq \frac{\langle v,\mathcal Lv\rangle}{(|v|^2+\tau^2)^{1/2}}
 \geq -C\tr a.
\]
Since $\rho_\tau=w=0$ at $t=0,1$, the maximum principle and
\eqref{sp:eq:strictbarrier} give $\rho_\tau\leq C'w$.
Letting $\tau\downarrow0$, we obtain
\begin{equation}\label{sp:eq:vbarrier}
 |d_Xu-d_Xh|\leq C'(u-\underline u).
\end{equation}
Dividing by $t$ or $1-t$
and let $t\uparrow 0$, we get
$|d_Xu_t-d_Xh_t|\leq C$ at both endpoints. Since $h$ is fixed,
this yields the desired bound for $|d_Xu_t|$ on the boundary.

We next bound $c=u_{tt}$ at the endpoints. At $t=0,1$, $A_u=A_j$ determined by $\varphi_j$ and the mixed derivative estimate gives
$|b|=|\partial_X u_t|\leq C_b$. Choose $\delta_0>0$ so that
$A_j-\delta_0 I\in\Gamma_k$ for $j=0,1$, and set
\[
 m_0=\min_{x\in X,\;j=0,1}f(A_j(x)-\delta_0 I)>0.
\]
If $c\geq C_b^2/\delta_0$, then
\[
 M=A_j-\frac{bb^*}{c}\geq A_j-\delta_0 I,
 \qquad f(M)\geq m_0.
\]
The equation $cf(M)=\varepsilon$ then gives $c\leq\varepsilon/m_0$.
Thus $u_{tt}$ is uniformly bounded at both endpoints.

Finally, differentiating the equation twice with respect to $t$
and using the concavity of $\F$, we obtain
\[
 \mathcal Lu_{tt}
 =-D^2\F(B_u)[(B_u)_t,(B_u)_t]\geq0.
\]
The maximum principle gives $u_{tt}\leq C$ on $X\times[0,1]$.
Since $W=1/u_{tt}$ on $\Omega$, the last assertion follows.
\end{proof}

\subsection{The Real Hessian Of $Z$}
\label{rh:sec:real}

We now estimate $\nabla_X^2Z$. On the component of $\partial\Omega$
corresponding to $t=j$, the Legendre derivative formula gives
\[
 \nabla_X^2Z
 =\nabla_X^2\varphi_j
  -\frac{d_Xu_t\otimes d_Xu_t}{u_{tt}}.
\]
By Lemma \ref{sp:lem:boundarySchur},
\[
 \left|\frac{d_Xu_t\otimes d_Xu_t}{u_{tt}}\right|=\frac{|d_X u_t|^2}{c}\leq C\frac{|b|^2}{c}=C\left(\tr A_j-\tr M\right)\leq C\tr A_j.
\]
Together with the preceding estimates, this gives
\begin{equation}\label{rh:eq:inputs}
 |Z|+|\partial_XZ|+|\partial_X\bar\partial_XZ|\leq C,
 \qquad |Z_s|\leq1,
 \qquad |\nabla_X^2Z|_{\partial\Omega}\leq C.
\end{equation}

\begin{proposition}\label{rh:prop:real}
Under~\eqref{rh:eq:inputs}, we have $|\nabla_X^2Z|\leq C$ on
$\Omega$, with $C$ independent of $\varepsilon$.
\end{proposition}

\begin{proof}
We follow the real Hessian estimate of Chu--McCleerey
\cite[Section 4]{CM} and treat the $s$ derivatives separately.
Set
\[
 N=\nabla_X^2Z,\qquad d=2n,\qquad
 L=1+\sup_{\overline\Omega}|N|.
\]
Since $\tr N\geq-C_0$, we have
\[
 \lambda_{\min}(N)
 \geq -C_0-(d-1)\lambda_{\max}(N).
\]
Hence
\[
 L\leq C\left(1+
 \sup_{\overline\Omega}\max\{\lambda_{\max}(N),0\}\right).
\]

Set $P=|\partial_XZ|^2$, $\rho=N+Lg$ and $R=|\rho|^2$. Then $\rho>0$ and
$R\leq4dL^2$. After subtracting a uniform constant from $Z$, we may assume
$-C_Z\leq Z\leq-1$. For a constant $A\geq2$ to be chosen below, on the set where $\lambda_{\max}(N)>0$, we consider
\begin{equation}\label{rh:eq:Q}
 Q=\log\Lambda+\xi(R)+\eta(P)+e^{-AZ},
 \qquad \Lambda=\lambda_{\max}(N),
\end{equation}
where \[
 \xi(R)=-\frac13\log\left(1-\frac{R}{5dL^2}\right),
 \qquad
 \eta(P)=-\frac13\log(1+\sup_\Omega P-P).
\]
It is easy to check that
\[
 \xi''=3(\xi')^2,\qquad \eta''=3(\eta')^2,\qquad
 \frac{c}{L^2}\leq\xi'\leq\frac{C}{L^2},
 \qquad c\leq\eta'\leq C.
\]

The sum $\xi(R)+\eta(P)+e^{-AZ}$ in $Q$ have bounded oscillation for fixed $A$.
Thus, if $Q$ attains its maximum on $\partial\Omega$,
\eqref{rh:eq:inputs} gives the desired estimate. We may assume
that $Q$ attains its maximum at an interior point $(x_0,s_0)$.
Comparing the values of $Q$ at $(x_0,s_0)$ and at a maximum
point of $\Lambda$, we obtain
\[
 L\leq C_A\bigl(1+\Lambda(x_0,s_0)\bigr).
\]
If $\Lambda(x_0,s_0)<1$, this already gives the desired bound.
We may therefore assume $\Lambda(x_0,s_0)\geq1$, in which case
\begin{equation}\label{rh:eq:Llambda}
 L\leq C_A\Lambda(x_0,s_0).
\end{equation}

\smallskip
\noindent\emph{Calculations at $(x_0,s_0)$.}

Choose a unit eigenvector $V_1$ of $N$ corresponding to $\Lambda$
at $(x_0,s_0)$. Extend it to a unit vector field $V_1(x)$ near $x_0$,
independent of $s$, with $\nabla_X V_1(x_0)=0$ and bounded second
covariant derivatives at $x_0$. Let $\widetilde\Lambda$ be the
largest eigenvalue of
\[
 \widetilde N
 =N-S
\]
where $S=\bigl(g-g(V_1,\cdot)\otimes g(V_1,\cdot)\bigr)$. At $(x_0,s_0)$, the perturbation leaves the eigenvalue in the
$V_1$ direction unchanged and decreases every eigenvalue on
$V_1^\perp$ by $1$. Hence $\widetilde\Lambda=\Lambda$ at this
point, and $\widetilde\Lambda$ is simple and smooth.
Since $\widetilde\Lambda\leq\Lambda$ nearby, the function $\widetilde Q$ obtained
by replacing $\log\Lambda$ with $\log\widetilde\Lambda$ in $Q$
still has a local maximum at $(x_0,s_0)$. 

For notational convenience, we write $Q$ and $\Lambda$ for
$\widetilde Q$ and $\widetilde\Lambda$, respectively, in the
rest of the proof. Let $\lambda_\alpha$, $\alpha>1$, denote the other eigenvalues
of $\widetilde N$ at $(x_0,s_0)$, with corresponding
orthonormal eigenvectors $V_\alpha$.
Choose a unitary frame in which $M$ is diagonal at this point.
Order its eigenvalues so that $m_1\geq\cdots\geq m_n$ and hence
$F_1\leq\cdots\leq F_n$. Differentiating~\eqref{sp:eq:ZPDE}
twice along $X$ gives
\[
 \D N_{ab}=\B(M_a,M_b)+\mathcal E_{ab},
 \qquad
 |\mathcal E|\leq C(1+|N|)\A.
\]
where $\displaystyle\mathcal E=\sum_i F_i\left[
     (\nabla_i\nabla_{\bar i}N)_{ab}
     -(\nabla_a\nabla_b N)_{i\bar i}\right]$ and commuting the derivatives gives
\[|\mathcal E|\leq C(1+|N|)\A.\]
It follows that
\begin{equation}\label{rh:eq:R}
\begin{split}
    \D R&=2\langle\rho,\mathcal DN\rangle
  +2\sum_iF_i|N_i|^2+2\varepsilon|N_s|^2\\
 &\geq 2\sum_i F_i|N_i|^2+2\varepsilon|N_s|^2-CL^2\A,
\end{split}
\end{equation}
where we dropped the nonnegative contraction
$2\rho_{ab}\B(M_a,M_b)$. Together with~\eqref{sp:eq:DP} and the uniform bound for $P$, we also have
\begin{equation} \label{rh:eq:P}
\begin{split}
        \D P
 &\geq \sum_{i,p}F_i\bigl(|Z_{ip}|^2+|Z_{i\bar p}|^2\bigr)
       +2\varepsilon|\partial_X Z_s|^2-C\A.
\end{split}
\end{equation}

At the maximum point, $\nabla_X S=0$ and $S_s=0$. For the simple eigenvalue $\Lambda$, we have
\begin{align*}
 \mathcal D\Lambda
 ={}&\mathcal B(M_{V_1},M_{V_1})
 +\mathcal E(V_1,V_1)-(\mathcal D S)(V_1,V_1)\\
 &+2\sum_{i,\alpha>1}F_i
   \frac{|N_{V_1V_\alpha;i}|^2}{\Lambda-\lambda_\alpha}
 +2\varepsilon\sum_{\alpha>1}
   \frac{|N_{V_1V_\alpha;s}|^2}{\Lambda-\lambda_\alpha}.
\end{align*}
Here $N_{V_1V_\alpha;i}=(\nabla_i N)(V_1,V_\alpha)$, and
similarly for the $s$ derivative. Since the second derivatives of
$S$ are bounded and $S$ is independent of $s$, we have
\[
 |\mathcal E(V_1,V_1)-(\mathcal D S)(V_1,V_1)|
 \leq C(1+|N|)\mathcal A.
\]
At the maximum point, the lower bound for $\tr N$ gives
$|N|\leq C(1+\Lambda)$. Since $\Lambda\geq1$, at the maximum point, we have
\[
\begin{split}
     \mathcal D\Lambda
 \geq{}&\mathcal B(M_{V_1},M_{V_1})
 +2\sum_{i,\alpha>1}F_i
   \frac{|N_{V_1V_\alpha;i}|^2}{\Lambda-\lambda_\alpha}\\
 &+2\varepsilon\sum_{\alpha>1}
   \frac{|N_{V_1V_\alpha;s}|^2}{\Lambda-\lambda_\alpha}
 -C\Lambda\mathcal A.
\end{split}
\]
The error in $\xi'\mathcal D R$ has the same bound since
$\xi'\leq C/L^2$. These constants are independent of $A$.

The first derivative conditions along $s$ are
\[
 \frac{\Lambda_s}{\Lambda}
 =-\xi'R_s-\eta'P_s+Ae^{-AZ}Z_s.
\]
Using the second identity and
$\xi''=3(\xi')^2$, $\eta''=3(\eta')^2$, we obtain
\begin{align}
 &\varepsilon\left[
 -\frac{\Lambda_s^2}{\Lambda^2}
 +\xi''R_s^2+\eta''P_s^2+A^2e^{-AZ}Z_s^2
 \right]\notag\\
 &\hspace{12mm}\geq
 -3\varepsilon A^2e^{-2AZ}Z_s^2
 \geq -C_A\varepsilon.
 \label{rh:eq:time}
\end{align}
Here we used $|Z_s|\leq1$ and $-C_Z\leq Z\leq-1$.

Set
\[
 E=\sum_{i,p}F_i\bigl(|Z_{ip}|^2+|Z_{i\bar p}|^2\bigr),
 \qquad
 B=\sum_i F_i\frac{|N_{V_1V_1;i}|^2}{\Lambda^2},
\]
and define
\begin{align*}
 G_1&=2\sum_{i,\alpha>1}
       \frac{F_i|N_{V_1V_\alpha;i}|^2}
            {\Lambda(\Lambda-\lambda_\alpha)},\\
 G_2&=\frac1\Lambda\sum_{i\ne q}
       \frac{F_i-F_q}{m_q-m_i}|(M_{V_1})_{i\bar q}|^2,\\
 G_3&=\frac1{C_A\Lambda^2}\sum_iF_i|N_i|^2.
\end{align*}
By~\eqref{rh:eq:Llambda},
the term $G_3$ is bounded above by $2\xi'\sum_iF_i|N_i|^2$
if $C_A$ is chosen large enough. Using
$\D Z=-\A$, \eqref{rh:eq:R}--\eqref{rh:eq:time}, at $(x_0,s_0)$, we have
\begin{align}
0\geq \mathcal DQ
={}&
 \frac{\mathcal D\Lambda}{\Lambda}-B
 +\xi'\mathcal DR+\eta'\mathcal DP
 -Ae^{-AZ}\mathcal DZ
 \notag\\
&+\xi''\sum_iF_i|R_i|^2
 +\eta''\sum_iF_i|P_i|^2
 +A^2e^{-AZ}\sum_iF_i|Z_i|^2
 \notag\\
&+\varepsilon\left(
 -\frac{\Lambda_s^2}{\Lambda^2}
 +\xi''R_s^2+\eta''P_s^2
 +A^2e^{-AZ}Z_s^2
 \right)
 \notag\\
\geq{}&
 G_1+G_2-C\mathcal A-B
 +2\xi'\sum_iF_i|N_i|^2
 +\eta'E
 -(C\xi'L^2+C\eta')\mathcal A
 \notag\\
&+\xi''\sum_iF_i|R_i|^2
 +\eta''\sum_iF_i|P_i|^2
 +A^2e^{-AZ}\sum_iF_i|Z_i|^2
 +Ae^{-AZ}\mathcal A
 \notag\\
&+\varepsilon\left(
 -\frac{\Lambda_s^2}{\Lambda^2}
 +\xi''R_s^2+\eta''P_s^2
 +A^2e^{-AZ}Z_s^2
 +2\xi'|N_s|^2
 +2\eta'|\partial_X Z_s|^2
 \right)
 \notag\\
\geq{}&
 G_1+G_2+G_3-B+cE
 +\xi''\sum_iF_i|R_i|^2
 +\eta''\sum_iF_i|P_i|^2
 \notag\\
&+A^2e^{-AZ}\sum_iF_i|Z_i|^2
 +Ae^{-AZ}\mathcal A
 -C_0\mathcal A-C_A\varepsilon.
 \label{rh:eq:ledger}
\end{align}
Here we used $\mathcal DZ=-\mathcal A$.
In the last inequality, \eqref{rh:eq:Llambda} and
$\xi'\geq c/L^2$ give
\[
 2\xi'\sum_iF_i|N_i|^2
 \geq \frac{1}{C_A\Lambda^2}\sum_iF_i|N_i|^2
 =G_3.
\]

\smallskip
\noindent\emph{Estimate of the negative term.}

We next estimate $\displaystyle B=\sum_i F_i\frac{|N_{V_1V_1;i}|^2}{\Lambda^2}$ in~\eqref{rh:eq:ledger}. At  $(x_0,s_0)$, $Q_i=0$ gives
\[\frac{\Lambda_i}{\Lambda}
 =-\xi'R_i-\eta'P_i+Ae^{-AZ}Z_i.\]
Since $\Lambda_i=N_{V_1V_1;i}$ the inequality
$|a+b+c|^2\leq3(|a|^2+|b|^2+|c|^2)$ implies
\[
 B\leq
 \xi''\sum_iF_i|R_i|^2+\eta''\sum_iF_i|P_i|^2
 +3A^2e^{-2AZ}\sum_iF_i|Z_i|^2.
\]
Substituting this into~\eqref{rh:eq:ledger} and using
$|\partial_X Z|\leq C$, $\A\geq1$, and $\varepsilon\leq1$,
we get
\begin{equation}\label{rh:eq:coarse}
 E\leq C_A\A.
\end{equation}

Set $\tau=A^{-2}e^{2AZ}\leq1/4$. Let $i_0$ be the largest index
such that $F_{i_0}<\tau F_{i_0+1}$, taking $i_0=0$ if there
is no such index, and put $I=\{i_0+1,\ldots,n\}$. Then
\[
 F_i\geq\tau^{n-1}F_n\geq\frac{\tau^{n-1}}n\A
 \quad(i\in I),
 \qquad
 F_q\leq\tau F_i
 \quad(i\in I,\ q\notin I).
\]

For $i\in I$, the bound $F_i\geq \tau^{n-1}\mathcal A/n$
and \eqref{rh:eq:coarse} give
\[
 \sum_p\bigl(|Z_{ip}|^2+|Z_{i\bar p}|^2\bigr)
 \leq \frac{E}{F_i}\leq C_A,
\]
where we used $\tau\geq A^{-2}e^{-2AC_Z}$.
Write
\[
 W_1=\frac{V_1-\sqrt{-1}JV_1}{\sqrt2}
     =\sum_q\nu_q e_q.
\]
The eigenvector equation $NV_1=\Lambda V_1$ then gives
$|\nu_q|\leq C_A/\Lambda$ for $q\in I$.

For $0<\gamma<1/2$, write $B=B_1+B_2+B_3$, where
\begin{align*}
 B_1 &= \sum_{i\notin I}F_i
          \frac{|N_{V_1V_1;i}|^2}{\Lambda^2},\\
 B_2 &= 2\gamma\sum_{i\in I}F_i
          \frac{|N_{V_1V_1;i}|^2}{\Lambda^2},\\
 B_3 &= (1-2\gamma)\sum_{i\in I}F_i
          \frac{|N_{V_1V_1;i}|^2}{\Lambda^2}.
\end{align*}
Arguing as in the derivation of \eqref{rh:eq:coarse}, we obtain
\begin{align}
 B_1+B_2\leq{}&
 \xi''\sum_iF_i|R_i|^2+\eta''\sum_iF_i|P_i|^2\notag\\
 &+6\gamma A^2e^{-2AZ}\sum_iF_i|Z_i|^2+C\A.
 \label{rh:eq:B12}
\end{align}

It remains to control $B_3$. Write
$JV_1=\displaystyle\sum_{\alpha>1}\mu_\alpha V_\alpha$. Using $V_1=\sqrt{2}\,\overline W_1-\sqrt{-1}\,JV_1$,
we obtain
\begin{align*}
 N_{V_1V_1;i}
 &=-\sqrt{-1}\sum_{\alpha>1}\mu_\alpha N_{V_1V_\alpha;i}
   +\sqrt{2}\sum_q\overline{\nu_q}\,N_{V_1\bar q;i}\\
 &=-\sqrt{-1}\sum_{\alpha>1}\mu_\alpha N_{V_1V_\alpha;i}
   +\sqrt{2}\sum_q\overline{\nu_q}(M_{V_1})_{i\bar q}
   -\sqrt{2}\,d_XZ\!\left(
       \operatorname{Rm}(e_i,V_1)\overline W_1
     \right)\\
     &=-\sqrt{-1}\sum_{\alpha>1}\mu_\alpha N_{V_1V_\alpha;i}
  +\sqrt2\sum_q\overline{\nu_q}(M_{V_1})_{i\bar q}
  +O(1)
\end{align*}
Split the sum over $q$ into
$q\notin I$ and $q\in I$. By Young's inequality and the bound
for $\nu_q$, when $q\in I$, we have
\begin{align*}
 B_3&\leq B_{31}+B_{32}+B_{33}
               +\frac{C\A}{\gamma\Lambda^2},\\
 B_{31}
 &=\frac{1-\gamma}{\Lambda^2}\sum_{i\in I}F_i
   \left|\sum_{\alpha>1}\mu_\alpha N_{V_1V_\alpha;i}\right|^2,\\
 B_{32}
 &=\frac{C}{\gamma\Lambda^2}
   \sum_{i\in I,\ q\notin I}F_i|(M_{V_1})_{i\bar q}|^2,\\
 B_{33}
 &=\frac{C_A}{\gamma\Lambda^4}
   \sum_{i\in I,\ q\in I}F_i|(M_{V_1})_{i\bar q}|^2.
\end{align*}

For $B_{31}$, Cauchy--Schwarz gives
\[
 \left|\sum_{\alpha>1}\mu_\alpha N_{V_1V_\alpha;i}\right|^2
 \leq
 \left(\sum_{\alpha>1}(\Lambda-\lambda_\alpha)\mu_\alpha^2\right)
 \left(\sum_{\alpha>1}
       \frac{|N_{V_1V_\alpha;i}|^2}{\Lambda-\lambda_\alpha}\right).
\]
Since $JV_1\perp V_1$ and $|JV_1|=1$, we have
$S(JV_1,JV_1)=1$. Thus the first factor is
\[
 \sum_{\alpha>1}(\Lambda-\lambda_\alpha)\mu_\alpha^2
 =\Lambda-(N-S)(JV_1,JV_1)
 =2\Lambda-2Z_{W_1\overline W_1}+1
 \leq 2\Lambda+C.
\]
It follows that
\[
 B_{31}
 \leq \frac{(1-\gamma)(2\Lambda+C)}{2\Lambda}\,G_1
 \leq G_1
 \qquad\text{if }\Lambda\geq C/\gamma,
\]
after increasing $C$ if necessary.

For $i\in I$ and $q\notin I$, we have
$F_i-F_q\geq F_i/2$ and $0<m_q-m_i\leq C$.
Therefore the corresponding coefficient in $G_2$ is at least
$cF_i/\Lambda$, and $B_{32}\leq G_2$ when
$\Lambda\geq C/\gamma$.

For $B_{33}$, commuting derivatives gives
\[
 \sum_{i\in I,\ q\in I}F_i|(M_{V_1})_{i\bar q}|^2
 \leq C\sum_iF_i|N_i|^2+C\A.
\]
Substituting this bound into the definition of $B_{33}$ gives
\[
 B_{33}
 \leq \frac{C_A}{\gamma\Lambda^4}\sum_iF_i|N_i|^2
      +\frac{C_A}{\gamma\Lambda^4}\A.
\]
Since $G_3=(C_A'\Lambda^2)^{-1}\displaystyle\sum_iF_i|N_i|^2$,
the first term is at most $G_3$ if
$\Lambda^2\geq C_A C_A'/\gamma$. The second term is at
most $\A/\gamma$ if $\Lambda^4\geq C_A$.
Both conditions hold when $\Lambda\geq C_A''/\gamma$
for a sufficiently large $C_A''$.

Combining this with the bounds for $B_{31}$ and $B_{32}$
gives, after increasing $C_A$,
\begin{equation}\label{rh:eq:B3}
 B_3\leq G_1+G_2+G_3+\frac{C_0}{\gamma}\A
 \qquad\text{if }\Lambda\geq C_A/\gamma,
\end{equation}
where $C_0$ is independent of $A$.

\smallskip
\noindent\emph{Choice of constants.}
Combining~\eqref{rh:eq:ledger}, \eqref{rh:eq:B12},
and~\eqref{rh:eq:B3}, and increasing $C_0$ if necessary, we obtain
\[
 0\geq cE+Ae^{-AZ}\A
 +\bigl(A^2e^{-AZ}-6\gamma A^2e^{-2AZ}\bigr)
       \sum_iF_i|Z_i|^2
 -\frac{C_0}{\gamma}\A-C_A\varepsilon
\]
whenever $\Lambda\geq C_A/\gamma$. Fix $A=6C_0+2$ and take
$\gamma=e^{AZ}/6$ at the maximum point. Then $A^2e^{-AZ}-6\gamma A^2e^{-2AZ}=0$ and 
\[
 0\geq 2e^{-AZ}\A-C_A\varepsilon
\]
for sufficiently small $\varepsilon$, since $Z\leq-1$ and
$\A\geq1$. This contradiction shows that $\Lambda<C_A/\gamma$
at the maximum point. Since $\gamma\geq e^{-AC_Z}/6$,
\eqref{rh:eq:Llambda} yields the desired bound for $N$.

\end{proof}

\section{The Mixed Derivative Estimate}\label{mx:sec:mixed}

We use the notation of Subsection~\ref{sp:sec:legendre}
and write $N$ as in Proposition~\ref{rh:prop:real}.
The inverse Legendre formulas give
\begin{equation}\label{mx:eq:inverse}
 \nabla^2_{X,t}u=
 \begin{pmatrix}N&0\\0&0\end{pmatrix}
 +\frac1c
 \binom{\nabla_Xu_t}{c}\binom{\nabla_Xu_t}{c}^{T}.
\end{equation}
By \eqref{sp:eq:Legendre-identities}, the trace of the second term
is $2|p|^2/W+1/W$. Proposition~\ref{rh:prop:real} bounds $N$, and
Proposition~\ref{sp:prop:boundarymixed} gives $W\geq w_0>0$.
It therefore suffices to bound $|p|^2/W$.

\begin{proposition}\label{mx:prop:mixed}
There exist $C>0$ and $\varepsilon_0>0$ such that every smooth
strictly admissible solution of \eqref{eq:root} with
$0<\varepsilon\leq\varepsilon_0$ satisfies
\[
 \frac{|p|^2}{W}\leq C \quad\text{on }\Omega,
 \qquad
 |\nabla^2_{X,t}u|\leq C \quad\text{on }X\times[0,1].
\]
\end{proposition}
\begin{proof}
Set
\[
 P=|p|^2,\qquad S=W+1,\qquad q=P/S,\qquad
 T=\tr M,\qquad H=T+q.
\]
By Lemma~\ref{sp:lem:boundarySchur}, $q$ is bounded on
$\partial\Omega$. We use \eqref{rh:eq:inputs} and
Proposition~\ref{rh:prop:real} throughout. In particular,
$0<T\leq C$.
Choose a fixed $L_0$ such that $\rho_N=N+L_0g\geq g$, and set
$R_N=|\rho_N|^2$. In particular, $R_N$ is uniformly bounded.

Let $\underline u$ be the strict subsolution from
Subsection~\ref{sbsec:C0}. Set
$\Phi(x,s)=(\underline u-u)(x,t(x,s))$.
Consider
\begin{equation}\label{mx:eq:test}
 Q=\log H+\alpha R_N+a\Phi-AZ.
\end{equation}
We take $\alpha>0$ small, then $a>0$ small, and finally
$A>0$ large. All three constants are independent of $\varepsilon$;
the conditions on their sizes are given below.
Since $\alpha R_N+a\Phi-AZ$ has bounded oscillation, it suffices
to bound $H$ at a maximum point of $Q$. The boundary bound for $q$ also bounds
$H$ on $\partial\Omega$. Therefore, we just consider an interior maximum
with $H\geq\max\{1,2\sup_\Omega T\}$, so that
$q\geq H/2$ and $1/2\leq\rho=q/H<1$.

At the maximum point, choose a unitary frame diagonalizing $M$
and put $F=\operatorname{diag}(F_i)$. We use $\D$, $\A$, and $\B$
as in Subsection~\ref{sp:sec:complex}.
For the product rule, let the index $\mu$ range over the real
orthonormal directions $X_i,JX_i$ and $\partial_s$, where
$e_i=(X_i-\sqrt{-1}JX_i)/\sqrt2$. Set
$a_{X_i}=a_{JX_i}=F_i/2$ and $a_s=\varepsilon$.
For a real function $v$, we write
\[
 |\nabla v|_\D^2
 =\sum_\mu a_\mu|v_\mu|^2
 =\sum_i F_i|v_i|^2+\varepsilon v_s^2.
\]

By \eqref{sp:eq:DZ}, at this point we have
\begin{align}
 0&\geq\D Q
 =\frac{\D H}{H}-|\nabla\log H|_\D^2
      +\alpha\D R_N+a\D\Phi+A\A,\label{mx:eq:DQ}\\
 \nabla\log H&=-\alpha\nabla R_N-a\nabla\Phi+A\nabla Z.
 \notag
\end{align}

\medskip
\noindent\emph{Step 1. Estimates for the terms in $\D Q$.}
Differentiating \eqref{sp:eq:ZPDE} and commuting covariant
derivatives gives
\[
 \D W=-\B(M_s,M_s),\qquad
 \D p_j=\B(M_j,M_s)+E_j,\qquad
 \D T=\sum_j\B(M_j,M_j^*)+E_T.
\]
The first identity involves no commutation. In the second identity,
the commutator is curvature times $p$, so $|E|\leq C\A|p|$.
For $T$, contracting the commutator in the equation for
$\D N_{ab}$ gives $|E_T|\leq C\A$, since $N$ and $\partial_XZ$
are bounded. We also have
\[
 \D q=\frac{2}{S}\sum_\mu a_\mu
     \left|p_\mu-\frac{pW_\mu}{S}\right|^2
 +\frac2S\operatorname{Re}\sum_j\overline{p_j}\D p_j
 -\frac{P}{S^2}\D W.
\]
Adding $\D T$ combines the $\B$ terms into
$\sum_j\B(K_j,K_j^*)$, where $K_j=M_j+(p_j/S)M_s$.
Hence, for $H\geq1$,
\begin{equation}\label{mx:eq:DH}
 \frac{\D H}{H}\geq
 \frac1H\sum_j\B(K_j,K_j^*)
 +\frac2{SH}\sum_\mu a_\mu
       \left|p_\mu-\frac{pW_\mu}{S}\right|^2-C_0\A.
\end{equation}
where we use $|E|\leq C\mathcal A|p|$ and $|E_T|\leq C\mathcal A$. Since $N$ is now bounded, we obtain
\begin{align}
 \D R_N&=2\langle\rho_N,\D N\rangle
   +2\sum_\mu a_\mu|N_\mu|^2\notag\\
 &=2\sum_\mu a_\mu|N_\mu|^2
   +2\sum_{a,b}(\rho_N)_{ab}\B(M_a,M_b)
   +2\langle\rho_N,\mathcal E\rangle\notag\\
 &\geq 2\sum_\mu a_\mu|N_\mu|^2
    +2\sum_{a,b}(\rho_N)_{ab}\B(M_a,M_b)-C_1\A\notag\\
 &\geq 2\sum_\mu a_\mu|N_\mu|^2
    +4\sum_j\B(M_j,M_j^*)-C_1\A.\label{mx:eq:DR}
\end{align}
For the last inequality, the real matrix $(\B(M_a,M_b))$ is
positive semidefinite and $\rho_N\geq g$. We also use the fact that its trace in the real
directions is $\displaystyle 2\sum_j\B(M_j,M_j^*)$.

Set $r=s-(\underline u)_t(x,t(x,s))$. All derivatives of $\underline u$ are taken in the original
variables $(x,t)$ and evaluated at $(x,t(x,s))$.
The bounds in \eqref{eq:low} give $|r|+|\Phi|\leq C$.
Differentiating \eqref{sp:eq:ZPDE} in $s$ gives $\D t=0$.
Using the derivatives of $t$ from Subsection~\ref{sp:sec:legendre},
we have
\begin{align}
 \Phi_i&=\underline u_i-Z_i+rp_i,
 &\Phi_s&=-rW,\notag\\
 |\nabla\Phi|_\D^2&\leq C_\Phi(\A+J_p+J_t),
 &J_p&=p^*Fp,\qquad J_t=\varepsilon W^2.\label{mx:eq:Phigrad}
\end{align}
Using the definition of $Z$ and $\D Z=-\mathcal A$, we also have
\[
 \D\Phi
 =\tr\!\bigl(F(A_{\underline u}-b_{\underline u} p^*-p b_{\underline u}^*+c_{\underline u} pp^*)\bigr)
   +\varepsilon c_{\underline u} W^2-2\varepsilon W.
\]
Take $\eta$ as in the proof of Proposition~\ref{sp:prop:boundarymixed} such that $B_{\underline u}-\eta I\in\C_k$ everywhere.
By Lemma~\ref{lem:cone-equivalence}, we have
\[
 c'=c_{\underline u}-\eta>0,\qquad
 M'=A_{\underline u}-\eta I-b_{\underline u} b_{\underline u}^*/c'\in\Gamma_k.
\]
Then
\[
 A_{\underline u}-b_{\underline u} p^*-p b_{\underline u}^*+c_{\underline u} pp^*
 =\eta I+\eta pp^*+M'+c'(p-b_{\underline u}/c')(p-b_{\underline u}/c')^*.
\]
The last two terms have nonnegative pairing with $F$. We also have\[\varepsilon c_{\underline u} W^2-2\varepsilon W
 \geq (\eta/2)\varepsilon W^2-2\varepsilon/\eta=(\eta/2)J_t-2\varepsilon/\eta\]
 since \(c_{\underline u}\geq\eta\).
Since $\A\geq1$ by \eqref{sp:eq:elementary}, for sufficiently
small $\varepsilon$ we obtain
\begin{equation}\label{mx:eq:Phi}
 \D\Phi\geq c_\Phi(\A+J_p+J_t).
\end{equation}

\medskip
\noindent\emph{Step 2. A bound for $J_p+J_t$.}
The first derivative equation for $Q$ gives
\[
 |\nabla\log H|_\D^2
 \leq 3\alpha^2|\nabla R_N|_\D^2
        +3a^2|\nabla\Phi|_\D^2+3A^2|\nabla Z|_\D^2.
\]
Here $|\nabla R_N|_\D^2\leq C\sum_\mu a_\mu|N_\mu|^2$, since
$\nabla R_N=2\langle\rho_N,\nabla N\rangle$. By
\eqref{mx:eq:Phigrad} and the bounds for $\partial_X Z$ and $Z_s$,
\[
 |\nabla\Phi|_\D^2\leq C_\Phi(\A+J_p+J_t),
 \qquad |\nabla Z|_\D^2\leq C\A
\]
for $\varepsilon\leq1$. Substituting these in \eqref{mx:eq:DQ} and
using \eqref{mx:eq:DH}, \eqref{mx:eq:DR}, and \eqref{mx:eq:Phi},
we obtain
\[
 0\geq (2\alpha-C\alpha^2)\sum_\mu a_\mu |N_\mu|^2
          +(ac_\Phi-3a^2C_\Phi)(J_p+J_t)-C_A\A.
\]
Choose $\alpha$ and $a$ sufficiently small so that both coefficients
are positive. Then
\begin{equation}\label{mx:eq:weighted}
 J_p+J_t\leq C_A\A.
\end{equation}

\medskip
\noindent\emph{Step 3. The terms with derivatives along $X$.}
For a vector $v$ and a matrix $J$, we write
\[
 |v|_F^2=v^*Fv,\qquad \|J\|_F^2=\sum_{i,j}F_i|J_{ji}|^2.
\]
Since $q=P/S\geq H/2>0$ at the maximum point, $P>0$. Set
\[
 e=p/\sqrt P=\partial_X Z/\sqrt P,\quad V=M_s/\sqrt P=V^*,\quad
 U_{ji}=Z_{s,ji}/\sqrt P,\quad y_i=W_i/S.
\]
We also write
\[
 b=Ve,\qquad a_i=e^*U_{\cdot i},\qquad \Pi=I-ee^*.
\]
Then we have 
\begin{align*}
 (\log q)_i
 &=\frac{P_i}{P}-\frac{W_i}{S}\\
 &=\frac{1}{P}\sum_j
   \left(\overline{p_j}\,Z_{s,ji}+p_j\,Z_{s,i\bar j}\right)
   -\frac{W_i}{S}\\
 &=\sum_j\overline{e_j}\,U_{ji}
   +\sum_j V_{ij}e_j-y_i\\
 &=a_i+b_i-y_i.
\end{align*}
Set $d_q=\partial_X\log q$, so that $(d_q)_i=a_i+b_i-y_i$.
Set $v_\mu=p_\mu-pW_\mu/S$. By the definitions of $U$, $V$, and $y_i$,
\[
 v_i=\sqrt P\,(U_{\cdot i}-e y_i),
 \qquad
 v_{\bar i}=\sqrt P\,(V_{\cdot i}-e\overline{y_i}),
\]
where we used $W_{\bar i}=\overline{W_i}$. Since
$|v_{X_i}|^2+|v_{JX_i}|^2=|v_i|^2+|v_{\bar i}|^2$,
the terms with $\mu=X_i,JX_i$ in \eqref{mx:eq:DH} are
\begin{align*}
 \frac{2}{SH}\sum_i\frac{F_i}{2}
       \bigl(|v_{X_i}|^2+|v_{JX_i}|^2\bigr)
 &=\frac{1}{SH}\sum_iF_i
       \bigl(|v_i|^2+|v_{\bar i}|^2\bigr)\\
 &=\frac{P}{SH}\sum_iF_i\bigl(
       |U_{\cdot i}-e y_i|^2
       +|V_{\cdot i}-e\overline{y_i}|^2\bigr)\\
 &=\rho\mathcal S,
\end{align*}
where $\rho=q/H=P/(SH)$ and $\displaystyle\mathcal S=\sum_iF_i\bigl(
       |U_{\cdot i}-e y_i|^2
       +|V_{\cdot i}-e\overline{y_i}|^2\bigr)$.

We claim that
\begin{equation}\label{mx:eq:projection}
 \mathcal S-|d_q|_F^2
 =\|V\|_F^2+|b-y|_F^2+\|\Pi U\|_F^2
                           -2\operatorname{Re}(d_q^*Fb).
\end{equation}
The projection of $U_{\cdot i}$ onto $e$ is $a_i$.
Hermitian symmetry gives $e^*V_{\cdot i}=\overline b_i$.
Consequently the sum of the two squares in $\mathcal S$ for each $i$ equals
\[
 |d_{q,i}-b_i|^2+|\Pi U_{\cdot i}|^2
 +|b_i-y_i|^2+|\Pi V_{\cdot i}|^2.
\]
Use $|b_i|^2+|\Pi V_{\cdot i}|^2=|V_{\cdot i}|^2$ and sum with
weights $F_i$ to prove \eqref{mx:eq:projection}. 

If $d_i=H_i/H$ and $\tau_i=T_i/H$, then $d=\partial_X T/H+\partial_X q/H=\tau+\rho d_q$ and \eqref{mx:eq:projection} implies
\begin{align}
 \rho\mathcal S-|d|_F^2={}&
 \rho\bigl(\|V\|_F^2+|b-y|_F^2+\|\Pi U\|_F^2\bigr)
 -2\operatorname{Re}(d^*Fb)+2\operatorname{Re}(\tau^*Fb)\notag\\
 &+\frac{1-\rho}{\rho}|d|_F^2
 -\frac2\rho\operatorname{Re}(d^*F\tau)+\frac1\rho|\tau|_F^2.
                                                        \label{mx:eq:Hprojection}
\end{align}

We now apply \eqref{mx:eq:Hprojection} at the maximum point. Set
\[
 E_X=\sum_iF_i|N_i|^2,\qquad
 \beta=(A+a)\partial_XZ-a\partial_X\underline u-ar p=:\beta_0-ar p.
\]
Then $d=\beta-\alpha\partial_X R_N$ by $\partial_X Q=0$, $|\beta_0|\le C_A$, and
\begin{equation}\label{mx:eq:tracegradient}
 |\partial_XR_N|_F^2\le C E_X,\qquad
 |\tau|_F^2\le C H^{-2}E_X.
\end{equation}
By Cauchy--Schwarz on each row of $V$ and Hermitian symmetry,
\[
 |b|_F^2=|Ve|_F^2
 \leq\sum_{i,j}F_i|V_{ij}|^2=\|V\|_F^2.
\]
After dropping the nonnegative squares in \eqref{mx:eq:Hprojection},
we estimate the remaining cross terms by Young's inequality:
\begin{align*}
 2\alpha|(\partial_X R_N)^*Fb|
 &\leq \tfrac1{16}\|V\|_F^2+C\alpha^2 E_X,\\
 2|\tau^*Fb|
 &\leq \tfrac1{16}\|V\|_F^2+CH^{-2}E_X,\\
 \frac{2}{\rho}|d^*F\tau|
 &\leq \left(\frac{\alpha}{4}+\frac{C\alpha}{H}\right)E_X
       +\frac{C}{\alpha H^2}|\beta|_F^2.
\end{align*}
For the last line, just split $d=\beta-\alpha\partial_X R_N$ and use
$\rho\geq1/2$. For $\alpha$ small and $H$ large, the coefficients
of $E_X$ sum to at most $\alpha$.
In particular,
\begin{equation}\label{mx:eq:spatiallower}
 \rho\mathcal S-|d|_F^2+2\alpha E_X
 \geq \tfrac14\|V\|_F^2-2\operatorname{Re}(\beta^*FVe)
                    -\frac C{\alpha H^2}|\beta|_F^2,
\end{equation}

We deal the terms containing $\mathcal B$ in \eqref{mx:eq:DH}. For each $j$, complete the square in $\B$ using
$K_j=M_j+(p_j/S)M_s$. This gives
\begin{align}
 \frac1H\sum_j\B(K_j,K_j^*)+4\alpha\sum_j\B(M_j,M_j^*)&=\left(\frac1H+4\alpha\right)\mathcal B(Y_j,Y_j^*)
+\frac{4\alpha}{1+4\alpha H}\frac P{S^2}\B(M_s,M_s)\notag\\
 &\ge\frac{4\alpha}{1+4\alpha H}\frac P{S^2}\B(M_s,M_s)\notag\\
 &=\frac{4\alpha q^2}{1+4\alpha H}\B(V,V)
 \ge\frac q4\B(V,V)\label{mx:eq:concavitybudget}
\end{align}
where $Y_j=M_j+\frac{(p_j/S)M_s}{1+4\alpha H}$. The final inequality uses $1\leq 4\alpha H$ and $q/H\ge1/2$.

Then we estimate the term $\beta^*FVe$ in \eqref{mx:eq:spatiallower}. Write $V=D+O$ in the eigenframe of
$M$, where $D$ is diagonal and $O$ is off-diagonal. Order the
eigenvalues of $M$ as $\lambda_1\geq\cdots\geq\lambda_n$. We first
show that
\begin{equation}\label{mx:eq:offdiag}
 \|O\|_F^2\leq\frac{n-1}{k-1}\lambda_1\B(V,V).
\end{equation}
Put $\zeta=k^{-1}\sigma_k^{1/k-1}$ and
$w_i=\sigma_{k-1}(\lambda|i)$. Then $F_i=\zeta w_i$.
For $i\ne j$, put
$b_{ij}=\sigma_{k-2}(\lambda|ij)>0$ and
$a_{ij}=\sigma_{k-1}(\lambda|ij)$. Expansion in the omitted entries gives
\[
 w_i+w_j=2a_{ij}+(\lambda_i+\lambda_j)b_{ij}.
\]
If $a_{ij}\leq0$, $w_i+w_j$ is at most $2\lambda_1b_{ij}$.
If $a_{ij}>0$, the vector $\mu=\lambda|ij$, of length $n-2$,
belongs to $\Gamma_{k-1}$. Deleting two entries from $\lambda\in\Gamma_k$ gives
$\sigma_\ell(\mu)>0$ for $1\leq\ell\leq k-2$.
Since $a_{ij}=\sigma_{k-1}(\mu)>0$, it follows that
$\mu\in\Gamma_{k-1}$.
The Newton inequalities imply
\[
 \frac{a_{ij}}{b_{ij}}=\frac{\sigma_{k-1}(\mu)}{\sigma_{k-2}(\mu)}
 \leq\frac{n-k}{k-1}\frac{\sigma_1(\mu)}{n-2}
 \leq\frac{n-k}{k-1}\lambda_1.
\]
Thus 
\[
 w_i+w_j\leq\frac{2(n-1)}{k-1}\lambda_1b_{ij}.
\]
Since $F_j-F_i=\zeta(\lambda_i-\lambda_j)b_{ij}$,
the off-diagonal terms in \eqref{eq:B-spec} sum to
$\displaystyle 2\zeta\sum_{i<j}b_{ij}|V_{ij}|^2$.
The diagonal part is nonnegative by concavity. Thus
\begin{align*}
 \|O\|_F^2
 &=\sum_{i<j}(F_i+F_j)|V_{ij}|^2\\
 &=\zeta\sum_{i<j}(w_i+w_j)|V_{ij}|^2\\
 &\leq \frac{n-1}{k-1}\lambda_1
        \left(2\zeta\sum_{i<j}b_{ij}|V_{ij}|^2\right)\\
 &\leq \frac{n-1}{k-1}\lambda_1\B(V,V).
\end{align*}
This proves \eqref{mx:eq:offdiag}.

Since $\lambda_1$ is bounded, \eqref{mx:eq:offdiag} gives
$\|O\|_F^2\leq C\B(V,V)$. For any vector $z$, we also have
\begin{align}
 |z^*FDe|^2&\le |z|^2(e^*Fe)\|D\|_F^2,\label{mx:eq:diagcross}\\
 |z^*FOe|^2&\le |z|^2\A\|O\|_F^2.\label{mx:eq:offcross}
\end{align}
Using \eqref{mx:eq:diagcross}, \eqref{mx:eq:offcross}, and
Young's inequality, we obtain
\begin{align*}
 2|\beta_0^*FDe|&\le\tfrac1{16}\|D\|_F^2+C_A e^*Fe,\\
 2|\beta_0^*FOe|&\le\tfrac q8\B(V,V)+C_A\A/q,\\
 2a|r|\,|p^*FVe|&\le\tfrac1{16}\|V\|_F^2+C a^2r^2 J_p.
\end{align*}
Combining \eqref{mx:eq:spatiallower} and \eqref{mx:eq:concavitybudget}
with the three bounds above for the terms in
$-2\operatorname{Re}(\beta^*FVe)$, and using
$e^*Fe=J_p/P\leq C_A\A/q$ and
$|\beta|_F^2\leq 2|\beta_0|^2\A+2a^2r^2J_p\leq C_A\A$
by \eqref{mx:eq:weighted} and $P=Sq\geq q$, we obtain
\begin{align*}
 &\rho\mathcal S-|d|_F^2+2\alpha E_X
 +\frac1H\sum_j\B(K_j,K_j^*)
 +4\alpha\sum_j\B(M_j,M_j^*)\\
 &\quad\geq
 \frac14\|V\|_F^2-\frac1{16}\|D\|_F^2-\frac1{16}\|V\|_F^2
 +\frac q8\B(V,V)\\
 &\qquad\quad
 -C_A e^*Fe-\frac{C_A\A}{q}
 -\frac{C}{\alpha H^2}|\beta|_F^2-Ca^2r^2J_p\\
 &\quad\geq -C_A\A(q^{-1}+H^{-2})-Ca^2r^2J_p,
\end{align*}
where we used $\|D\|_F^2\leq\|V\|_F^2$ and $\B(V,V)\geq0$
in the last line.

\medskip
\noindent\emph{Step 4. The terms with $s$ derivatives and the conclusion.}
In the $s$ direction, $Q_s=0$ gives
\[
 d_s=H_s/H=-\alpha(R_N)_s+A Z_s+arW.
\]
Consequently
\[
 -\varepsilon d_s^2\ge
 -C\alpha^2\varepsilon|N_s|^2-3A^2\varepsilon
                                  -3a^2r^2\varepsilon W^2.
\]
The first term on the right is absorbed by
$2\alpha\varepsilon|N_s|^2$ from \eqref{mx:eq:DR}.
The other term with an $s$ derivative in \eqref{mx:eq:DH} is nonnegative.

Choose $a$ small enough that
$Ca^2r^2J_p+3a^2r^2J_t\leq (ac_\Phi/2)(J_p+J_t)$.
Finally, take $A>C_0+\alpha C_1+3$. Combining the estimates in \eqref{mx:eq:DQ} gives
\begin{equation}\label{mx:eq:final}
 0\geq\D Q\geq (A-C_0-\alpha C_1)\A
 +\tfrac12a c_\Phi(J_p+J_t)
 -C_A\A(q^{-1}+H^{-2})-3A^2\varepsilon.
\end{equation}
If $\varepsilon\leq(6A^2)^{-1}$ and $H$ is sufficiently large,
then $3A^2\varepsilon\leq1/2$ and
$C_A(q^{-1}+H^{-2})\leq1$, since $q\geq H/2$.
As $A-C_0-\alpha C_1>3$ and $\A\geq1$ by
\eqref{sp:eq:elementary}, this gives
$\D Q>0$, a contradiction. This bounds $H$ at the maximum point of $Q$,
and hence throughout $\Omega$.

Since $W\geq w_0$, we have
$P/W=(1+1/W)q\leq (1+w_0^{-1})q$.
The full Hessian bound follows from \eqref{mx:eq:inverse}.
\end{proof}

\section{Existence Of The Regularized Solutions}\label{sec:existence}
We have proved the following a priori estimate.
\begin{proposition}\label{prop:apriori}
There are constants $C$ and $\varepsilon_0>0$, depending only on
the background and $\varphi_0,\varphi_1$, such that every smooth admissible
solution of~\eqref{eq:root}, with $0<\varepsilon\leq\varepsilon_0$,
satisfies
\begin{equation}\label{eq:apriori}
 \|u\|_{C^1}+\|\nabla^2_{X,t}u\|_\infty\leq C,
 \qquad |M|\leq C,\qquad 0<c\leq C.
\end{equation}
The same constants can be chosen for endpoints
$\tau\varphi_0,\tau\varphi_1$, with $0\leq\tau\leq1$.
\end{proposition}
\begin{proof}
The first assertion follows from the estimates in
the previous sections and the inverse
Legendre formula at the beginning of Section~\ref{mx:sec:mixed}.
For the second assertion, let
$A_j=I+(\varphi_j)_{i\bar j}$. Choose $\delta\in(0,1/4)$ so that
$A_j-2\delta I\in\Gamma_k$ for $j=0,1$. Convexity of $\Gamma_k$ gives
\[
 (1-\tau)I+\tau A_j-2\delta I
  =(1-\tau)(1-2\delta)I+\tau(A_j-2\delta I)\in\Gamma_k.
\]
Thus the endpoint matrices $A_{\tau\varphi_j}$ have a uniform strict margin. Their smooth
norms are also uniformly bounded. The constants $K$ and $\eta$ in the barrier construction and
the proof of Proposition~\ref{sp:prop:boundarymixed} can be
chosen independently of $\tau$. The maximum principle estimates are therefore uniform in $\tau$,
and the gradient estimate is also uniform because the endpoint data
range over a compact family and their rescaled values converge to
constants along a subsequence.
\end{proof}

\begin{lemma}\label{lem:fixedellipticity}
For fixed $\varepsilon\in(0,\varepsilon_0]$, the augmented blocks $B_u$ of
all solutions in Proposition~\ref{prop:apriori} lie in a compact
subset of $\C_k$. The equation is uniformly elliptic near that subset.
Its ellipticity constants may depend on $\varepsilon$.
\end{lemma}
\begin{proof}
By $|M|\leq C$, one has $f(M)\leq C_f$. The equation gives
\[
 c\geq\varepsilon/C_f,\qquad f(M)=\varepsilon/c\geq\varepsilon/C.
\]
The Maclaurin inequalities give a positive lower bound for
$\sigma_j(M)$, for every $1\leq j\leq k$, from this lower bound for
$\sigma_k(M)$. Together with $|M|\leq C$, these inequalities keep
$M$ in a compact subset of $\Gamma_k$.
The full Hessian bound controls $A,b,c$. The Schur formula and the
lower bound for $c$ now keep $B_u$ in a compact subset of $\C_k$.
Lemma~\ref{lem:rootoperator} and compactness give the ellipticity
claim.
\end{proof}

\begin{proposition}\label{prop:existence}
For each $0<\varepsilon\leq\varepsilon_0$, equation~\eqref{eq:root}
has a unique smooth admissible solution.
\end{proposition}

\begin{proof}
Fix $\varepsilon\in(0,\varepsilon_0]$ and set
\[
 \mathcal I_\varepsilon
 =\left\{\tau\in[0,1]\;\middle|\;
 \begin{aligned}
 &\text{there exists a smooth strictly admissible }u\\
 &\text{such that }\F(B_u)=\sqrt{\varepsilon},\\
 &\text{and }u|_{t=j}=\tau\varphi_j\quad (j=0,1)
 \end{aligned}
 \right\}.
\]
We show that $\mathcal I_\varepsilon=[0,1]$.

\smallskip
\noindent\emph{Nonemptiness.}
For $\tau=0$, set
\[
 u^{(0)}(x,t)=\frac{\varepsilon}{2f(I)}t(t-1).
\]
Its boundary values are zero and
$B_{u^{(0)}}=\operatorname{diag}(I,\varepsilon/f(I))\in\C_k$.
Thus $\F(B_{u^{(0)}})=\sqrt{\varepsilon}$, and
$0\in\mathcal I_\varepsilon$.

\smallskip
\noindent\emph{Openness.}
Let $\tau_0\in\mathcal I_\varepsilon$ and let $u$ be a
corresponding solution. The linearization at $u$ is
\[
 \mathcal L_u v
 =D\F(B_u)\left[
 \begin{pmatrix}
   (v_{i\bar j}) & (v_{it})\\
   (v_{t\bar j}) & v_{tt}
 \end{pmatrix}
 \right].
\]
This operator is uniformly elliptic
for the fixed solution $u$. By \cite{GT}, it is an isomorphism from the subspace
of $C^{2,\alpha}(X\times[0,1])$ with zero boundary values onto
$C^\alpha(X\times[0,1])$.
The implicit function
theorem yields strictly admissible solutions for $\tau$ near
$\tau_0$. Schauder regularity makes these solutions smooth.
Hence $\mathcal I_\varepsilon$ is relatively open in $[0,1]$.

\smallskip
\noindent\emph{Closedness.}
Let $\tau_m\in\mathcal I_\varepsilon$ converge to $\tau_\infty$,
and let $u_m$ be corresponding solutions.
Proposition~\ref{prop:apriori} gives a uniform $C^2$ bound, and
Lemma~\ref{lem:fixedellipticity} places the matrices $B_{u_m}$ in
a fixed compact subset $\mathcal K_\varepsilon\Subset\C_k$.
On this set, $\F$ agrees with the globally concave and uniformly
elliptic extension
\[
 \widetilde\F(B)
 =\inf_{B_0\in\mathcal K_\varepsilon}
 \bigl\{\F(B_0)+D\F(B_0)[B-B_0]\bigr\}.
\]

The interior Evans--Krylov estimates \cite{Evans,Krylov} and the
boundary estimates of Silvestre--Sirakov
\cite{SS} give a uniform $C^{2,\alpha}$
bound for some $\alpha\in(0,1)$. Schauder estimates then give
higher bounds up to the boundary. After passing to a subsequence,
$u_m$ converges smoothly to a solution with boundary values
$\tau_\infty\varphi_j$. Its matrices lie in $\mathcal K_\varepsilon$,
so the limit is strictly admissible. Thus
$\tau_\infty\in\mathcal I_\varepsilon$.

Since $\mathcal I_\varepsilon$ is nonempty, relatively open, and
relatively closed in $[0,1]$, it equals $[0,1]$.

Uniqueness follows by linearizing along the segment between two
solutions, which lies in $\C_k$ by Lemma~\ref{lem:cone-convexity},
and applying the maximum principle to their difference.
\end{proof}

\section{Comparison and the Weak Limit}\label{sec:limit}

Let $\C_k^*$ be the dual cone of $\overline{\C_k}$ with
respect to the trace pairing, and set
\[
 \mathcal Q_k=\{Q\in\C_k^*:\tr Q=1\}.
\]
Since $\C_k$ contains every positive definite matrix,
each $Q\in\C_k^*$ is positive semidefinite. Thus $\mathcal Q_k$
is compact.

For a Hermitian matrix $B$, we define
\[
 F_0(B)=\min_{Q\in\mathcal Q_k}\tr(QB).
\]
By duality, $F_0(B)\geq0$ if and only if $B\in\overline{\C_k}$.
Moreover, $F_0(B+sI)=F_0(B)+s$ for every $s\in\mathbb R$.
It follows that $F_0(B)>0$ precisely on $\C_k$, and
$F_0(B)=0$ precisely on $\partial\C_k$.

We first prove a comparison principle for
\begin{equation}\label{eq:bellman}
 F_0(B_u)=0.
\end{equation}

\begin{definition}\label{def:subgeodesic}
A bounded upper semicontinuous function $v$ on $X\times(0,1)$
is an admissible subgeodesic if every smooth upper test $\phi$
for $v$ satisfies $B_\phi\in\overline{\C_k}$ at the point of
contact. Equivalently, $v$ is a viscosity subsolution of
\eqref{eq:bellman}. Similarly, we can define an admissible supergeodesic.

Let $\mathcal S(\varphi_0,\varphi_1)$ be the set of all such $v$
satisfying
\[
 \limsup_{(y,t)\to(x,j)}v(y,t)\leq\varphi_j(x),
 \qquad x\in X,\quad j=0,1.
\]
We define the envelope on $X\times(0,1)$ by
\[
 U_{\mathrm{env}}(x,t)
 =\sup_{v\in\mathcal S(\varphi_0,\varphi_1)}v(x,t).
\]
\end{definition}

\begin{lemma}\label{lem:comparison}
Let $v$ be an admissible subgeodesic of
\eqref{eq:bellman} and $w$ an admissible suppergeodesic. Suppose that, for $j=0,1$,
\[
 \limsup_{(y,t)\to(x,j)}v(y,t)
 \leq \varphi_j(x)
 \leq \liminf_{(y,t)\to(x,j)}w(y,t).
\]
Then $v\leq w$ on $X\times(0,1)$.
\end{lemma}

\begin{proof}
Take $\eta>0$ as in the proof of
Proposition~\ref{sp:prop:boundarymixed}. Since
$B_{\underline u}-\eta I\in\C_k$, the definition of $F_0$ gives
$F_0(B_{\underline u})\geq\eta$. For $0<\theta<1$, set
$v_\theta=(1-\theta)v+\theta\underline u$. If $\phi$ is an upper
test for $v_\theta$, then
$\psi=(\phi-\theta\underline u)/(1-\theta)$ is an upper test for $v$.
By concavity,
\begin{equation}\label{eq:strictcomparison}
 F_0(B_\phi)
 \geq (1-\theta)F_0(B_\psi)+\theta F_0(B_{\underline u})
 \geq\theta\eta.
\end{equation}
Also, $v_\theta\leq w$ in the boundary limit.

In a local unitary frame, we write
\[
 \tr(QB_\phi)=\tr(a_Q(z)D^2\phi)+r_Q,
 \qquad r_Q=\tr Q_{XX},
\]
where $z$ denotes the real product coordinates on $X\times [0,1]$. If $E(z)$ is the
matrix converting coordinate directions to the chosen frame, then
$a_Q=\operatorname{Re}(E^*QE)$. Taking the real and imaginary parts
of $Q^{1/2}E$ gives $a_Q=\Sigma_Q\Sigma_Q^T$, with
\begin{equation}\label{eq:diffusion}
 |\Sigma_Q|+|D\Sigma_Q|\leq C
 \qquad\text{uniformly for }Q\in\mathcal Q_k
\end{equation}
on each smaller chart, since $\tr Q=1$. Thus $-F_0(B_\phi)$
satisfies the local structure condition of
\cite[Example 3.6]{CIL}.

Suppose that $v_\theta-w$ attains a positive maximum
$m$ at an interior point $z_0$. Choose a coordinate ball
$B_R(z_0)$ with closure in the chart and set
\[
 \widetilde v(z)
 =v_\theta(z)-m-\kappa\left(|z-z_0|^2-\frac{R^2}{2}\right).
\]
For sufficiently small $\kappa>0$, $\widetilde v$ remains a
strict subsolution: adding a constant does not change the
equation, and the quadratic term reduces the margin in
\eqref{eq:strictcomparison} by at most $C\kappa$.
Since $v_\theta-w\leq m$, we have
\[
 \widetilde v-w\leq-\frac{\kappa R^2}{2}
 \quad\text{on }\partial B_R(z_0),
 \qquad
 (\widetilde v-w)(z_0)=\frac{\kappa R^2}{2}>0.
\]
This contradicts strict comparison on $B_R(z_0)$
\cite[Section 5.C]{CIL}. Thus $v_\theta\leq w$, and letting
$\theta\downarrow0$ gives $v\leq w$.
\end{proof}

Let $u_\varepsilon$ be the solution of \eqref{eq:root} given by
Proposition~\ref{prop:existence}. To pass to the limit in
\eqref{eq:bellman}, we first show that
\begin{equation}\label{eq:residual}
 0<F_0(B_{u_\varepsilon})
 \leq f(I)^{-1/2}\sqrt{\varepsilon}
 =:r_\varepsilon.
\end{equation}
Strict admissibility gives $F_0(B_{u_\varepsilon})>0$.
For any $0<d<F_0(B_{u_\varepsilon})$, we have
$B_{u_\varepsilon}-dI\in\C_k$. By concavity and homogeneity of
$\F$,
\[
 \sqrt{\varepsilon}
 =\F(B_{u_\varepsilon})
 \geq \F(B_{u_\varepsilon}-dI)+\F(dI)
 >d\sqrt{f(I)}.
\]
Letting $d\uparrow F_0(B_{u_\varepsilon})$ proves
\eqref{eq:residual}.

\begin{proof}[Proof of Theorem~\ref{thm:main}]
By \eqref{eq:apriori}, a subsequence of $u_\varepsilon$ converges
in $C^{1,\beta}(X\times[0,1])$, for every $0<\beta<1$, to a function
$U\in C^{1,1}(X\times[0,1])$. The estimate in
\eqref{eq:apriori} passes to $U$, and $U|_{t=j}=\varphi_j$ for
$j=0,1$. By \eqref{eq:residual} and viscosity stability
\cite[Section 6]{CIL}, $U$ solves \eqref{eq:bellman}.
Lemma~\ref{lem:comparison} gives uniqueness for this Dirichlet
problem. Thus every subsequential limit is $U$, and the whole
family converges to $U$.

To identify the envelope, let $v$ be an admissible subgeodesic
with the prescribed boundary upper bounds. Its upper tests
satisfy $B_\phi\in\overline{\C_k}$, or equivalently
$F_0(B_\phi)\geq0$. Since $U$ is a supersolution of
\eqref{eq:bellman}, Lemma~\ref{lem:comparison} gives $v\leq U$.
On the other hand, $U$ is itself an admissible subgeodesic
with the prescribed boundary values. Hence $U$ is the envelope.

We next verify \eqref{eq:geodesic} almost everywhere. At each
point where $U$ has a second-order Taylor expansion, adding
and subtracting a small positive quadratic term gives upper
and lower tests. Letting the coefficient of that term tend
to zero gives $F_0(B_U)=0$, and hence $B_U\in\partial\C_k$.

At such a point, write $A=A_U$, $b=b_U$, and $c=U_{tt}$.
If $c>0$, the Schur complement $M=A-bb^*/c$ lies in
$\partial\Gamma_k$. Thus $\sigma_k(M)=0$, and
\eqref{eq:schurpoly} gives $P(B_U)=0$.
If $c=0$, the inequality $|b|^2\leq c\tr A$ on
$\overline{\C_k}$ gives $b=0$. The definition of $P$ then gives
$P(B_U)=0$ also in this case. Therefore $U$ satisfies
\eqref{eq:geodesic} almost everywhere.

Finally, for each fixed $t$, the functions $u_\varepsilon(\cdot,t)$
are strictly $k$-positive by admissibility. Since they converge
uniformly to $U(\cdot,t)$, the dual inequalities for
$\overline{\Gamma_k}$ pass to the limit in distributions.
Thus each slice $U(\cdot,t)$ is weakly $k$-positive.
\end{proof}

\section{The Path Distance}\label{sec:hessian-distance}

In this section, we prove Theorem~\ref{thm:hessian-distance} following the endpoint variation argument of Chen
\cite[Section 5]{Chen}.

\begin{proof}[Proof of Theorem~\ref{thm:hessian-distance}]
Write $\mu_v=\omega_v^k\wedge\omega^{n-k}$ and
$V=\int_X\omega^n$. Let $\gamma$ be a piecewise smooth path
in $\mathcal H_k$ joining $\varphi_0$ to $\varphi_1$.
For each $s$, let $u=u_\varepsilon(s,x,t)$ solve
\eqref{eq:root} with endpoints $\gamma_s$ and $\varphi_1$.
The endpoint pairs lie in a compact family, so the estimates
and the range of sufficiently small $\varepsilon$ can be
chosen uniformly, as in Proposition~\ref{prop:apriori}.
For fixed $\varepsilon$, the implicit function theorem
gives smooth dependence on $s$ on each smooth piece of $\gamma$.
Differentiating \eqref{eq:root} in $s$ and applying the maximum
principle gives
\[
 \|u_s(s,\cdot,\cdot)\|_\infty\leq\|\dot\gamma_s\|_\infty.
\]
The derivative $u_t$ is uniformly bounded by \eqref{eq:low}.

Let $\mathcal R_u$ be the top-degree form on the left-hand side of
\eqref{eq:geodesic}, evaluated at $u$. By \eqref{eq:schurpoly}
and $cf(M)=\varepsilon$, we have
\[
 \mathcal R_u
 =\binom nk^{-1}c\,\sigma_k(M)\,\omega^n
 =\binom nk^{-1}\varepsilon f(M)^{k-1}\omega^n.
\]
Since $M\leq A_u$ and $\int_X\mu_u=V$, H\"older's inequality
gives $\int_X\mathcal R_u\leq C\varepsilon$.

Set
\[
 e_\varepsilon(s,t)=\int_X u_t^2\,\mu_u,\qquad
 E_\varepsilon(s)=\int_0^1 e_\varepsilon(s,t)\,dt.
\]
Differentiating and integrating by parts gives
\begin{align*}
 \partial_t e_\varepsilon(s,t)
 =2\int_X u_t\,\mathcal R_u.
\end{align*}
Writing $T_u=k\,\omega_u^{k-1}\wedge\omega^{n-k}$, we have
\begin{align*}
 E_\varepsilon'(s)
 &=\int_0^1\!\!\int_X
   \left(2u_tu_{st}\,\mu_u
   +u_t^2\sqrt{-1}\partial_X\bar\partial_Xu_s\wedge T_u\right)dt\\
 &=\int_0^1\!\!\int_X
   \left(2u_tu_{st}\,\mu_u
   +u_s\sqrt{-1}\partial_X\bar\partial_X(u_t^2)\wedge T_u\right)dt\\
 &=2\int_0^1\partial_t\left(\int_X u_su_t\,\mu_u\right)dt
   -2\int_0^1\!\!\int_X u_s
     \left(u_{tt}\mu_u-\sqrt{-1}\partial_Xu_t
                  \wedge\bar\partial_Xu_t\wedge T_u\right)dt\\
 &=2\left[\int_X u_su_t\,\mu_u\right]_{t=0}^{t=1}
   -2\int_0^1\!\!\int_X u_s\,\mathcal R_u\,dt.
\end{align*}
Using the bounds for $u_t$, $u_s$ and
$\int_X\mathcal R_u$, together with the boundary values of $u_s$,
we obtain
\begin{align*}
 |e_\varepsilon(s,t)-E_\varepsilon(s)|
 &\leq C\varepsilon,\\
 |E_\varepsilon'(s)|
 &\leq 2\|\dot\gamma_s\|_{\gamma_s}
         \sqrt{e_\varepsilon(s,0)}+C\varepsilon,
\end{align*}
where $\|\cdot\|_{\gamma_s}$ is the norm in \eqref{eq:metric}.
Consequently, for a fixed sufficiently large $C_0$,
\begin{equation}\label{eq:derivativeE}
    \left|\frac{d}{ds}\sqrt{E_\varepsilon(s)+C_0\varepsilon}\right|
 \leq \|\dot\gamma_s\|_{\gamma_s}+C\sqrt{\varepsilon}.
\end{equation}

At $s=1$, both endpoints are $\varphi_1$. Let
$m=\min_X f(A_{\varphi_1})>0$. The function
$\varphi_1+\varepsilon t(t-1)/(2m)$ is a subsolution of
\eqref{eq:root}. Comparison and convexity in $t$ give
\[
 \varphi_1+\frac{\varepsilon}{2m}t(t-1)
 \leq u_\varepsilon(1,\cdot,t)\leq\varphi_1.
\]
Taking the one-sided derivatives at $t=0,1$ and using
convexity once more yields $|u_t(1,\cdot,t)|\leq\varepsilon/(2m)$.
Thus $E_\varepsilon(1)\leq V\varepsilon^2/(4m^2)$.
Integrating the differential inequality \eqref{eq:derivativeE} in $s$ now gives
\[
 \limsup_{\varepsilon\downarrow0}\sqrt{E_\varepsilon(0)}
 \leq L_2(\gamma).
\]
On the other hand, $u_\varepsilon(0,\cdot,t)$ is a smooth path
joining $\varphi_0$ to $\varphi_1$. By Cauchy--Schwarz in $t$,
\[
 d_2(\varphi_0,\varphi_1)
 \leq L_2(u_\varepsilon(0,\cdot,\cdot))
 \leq\sqrt{E_\varepsilon(0)}.
\]
Taking the infimum over $\gamma$, we conclude that
$\sqrt{E_\varepsilon(0)}\to d_2(\varphi_0,\varphi_1)$.

For each fixed $t$, uniform convergence and the uniform Hessian
bound give weak convergence
$\mu_{u_\varepsilon(0,\cdot,t)}\rightharpoonup\mu_{U(\cdot,t)}$.
Since $u_t(0,\cdot,t)\to U_t(\cdot,t)$ uniformly, we obtain
\[
 \int_X |U_t(\cdot,t)|^2\,\mu_{U(\cdot,t)}
 =\lim_{\varepsilon\downarrow0}e_\varepsilon(0,t)
 =\lim_{\varepsilon\downarrow0}E_\varepsilon(0)
 =d_2(\varphi_0,\varphi_1)^2,
\]
where the second equality follows from
$|e_\varepsilon(0,t)-E_\varepsilon(0)|\leq C\varepsilon$.
This proves \eqref{eq:hessian-distance}

The line segment in $\mathcal H_k$ gives finiteness, and
reversal and concatenation of paths give symmetry and the
triangle inequality. Finally, convexity in $t$ gives
\[
 U_t(\cdot,0)\leq\varphi_1-\varphi_0\leq U_t(\cdot,1).
\]
Evaluating \eqref{eq:hessian-distance} at the two endpoints yields
\[
 d_2(\varphi_0,\varphi_1)^2\geq
 \max\left\{
 \int_X(\varphi_0-\varphi_1)_+^2\,\mu_{\varphi_0},\
 \int_X(\varphi_1-\varphi_0)_+^2\,\mu_{\varphi_1}
 \right\}.
\]
The endpoint measures are strictly positive, so the right-hand
side is positive whenever $\varphi_0\ne\varphi_1$. Hence, $(\mathcal H_k,d_2)$ is a metric space.
\end{proof}

\section*{acknowledge}

The author first considered this problem seven years ago and
thanks Wei Wei for discussions at that time. The main difficulties
then were the gradient and complex Hessian estimates. 

When the author recently returned to the problem, chatGPT suggested
using a partial Legendre transform. GPT also assisted in  language editing and review.
The author takes full responsibility for the contents of the
paper, including all mathematical statements and proofs.

While this manuscript is completed, the author became aware of the preprint
\emph{Geodesics for generalised complex Hessian equations} by
Kartick Ghosh--ZiYue Wang--YuLun Xu \cite{GWX}
(\href{https://arxiv.org/abs/2611.25768}{arXiv:2611.25768}).
Their work overlaps with ours regarding the existence of weak geodesics and the induced path distance, and addresses a more general class of Hessian equations. For the $k$-Hessian equation
considered here, we prove a uniform bound for all real second
derivatives of the regularized solutions, which gives full
$C^{1,1}$ regularity of the weak geodesic.


\begin{thebibliography}{23}
\providecommand{\natexlab}[1]{#1}
\providecommand{\url}[1]{\texttt{#1}}
\expandafter\ifx\csname urlstyle\endcsname\relax
  \providecommand{\doi}[1]{doi: #1}\else
  \providecommand{\doi}{doi: \begingroup \urlstyle{rm}\Url}\fi

\bibitem[{\AA}hag et~al.(2024){\AA}hag, Czy{\.{z}}, Lu, and Rashkovskii]{ACLR}
Per {\AA}hag, Rafa{\l} Czy{\.{z}}, Chinh~H. Lu, and Alexander Rashkovskii.
\newblock Kiselman minimum principle and rooftop envelopes in complex {Hessian}
  equations.
\newblock \emph{Math. Z.}, 308:\penalty0 Paper No. 70, 2024.
\newblock \doi{10.1007/s00209-024-03624-3}.


\bibitem[Caffarelli et~al.(1985)Caffarelli, Nirenberg, and Spruck]{CNS}
Luis Caffarelli, Louis Nirenberg, and Joel Spruck.
\newblock The {Dirichlet} problem for nonlinear second order elliptic
  equations, {III}: Functions of the eigenvalues of the {Hessian}.
\newblock \emph{Acta Math.}, 155:\penalty0 261--301, 1985.
\newblock \doi{10.1007/BF02392544}.

\bibitem[Chen(2000)]{Chen}
Xiuxiong Chen.
\newblock The space of {K\"ahler} metrics.
\newblock \emph{J. Differential Geom.}, 56\penalty0 (2):\penalty0 189--234,
  2000.
\newblock \doi{10.4310/jdg/1090347643}.

\bibitem[Chen and He(2011)]{ChenHe}
Xiuxiong Chen and Weiyong He.
\newblock The space of volume forms.
\newblock \emph{Int. Math. Res. Not. IMRN}, 2011\penalty0 (5):\penalty0
  967--1009, 2011.

\bibitem[Chen and He(2019)]{ChenHeClass}
Xiuxiong Chen and Weiyong He.
\newblock A class of fully nonlinear equations.
\newblock \emph{Pure Appl. Math. Q.}, 15\penalty0 (4):\penalty0 1029--1045,
  2019.

\bibitem[Chu(2019)]{ChuVolume}
Jianchun Chu.
\newblock {$C^{1,1}$} regularity of geodesics in the space of volume forms.
\newblock \emph{Calc. Var. Partial Differential Equations}, 58\penalty0
  (6):\penalty0 Paper No. 194, 2019.
\newblock \doi{10.1007/s00526-019-1648-3}.

\bibitem[Chu and McCleerey(2021)]{CM}
Jianchun Chu and Nicholas McCleerey.
\newblock Fully non-linear degenerate elliptic equations in complex geometry.
\newblock \emph{J. Funct. Anal.}, 281\penalty0 (9):\penalty0 Paper No. 109176,
  2021.
\newblock \doi{10.1016/j.jfa.2021.109176}.

\bibitem[Chu et~al.(2017)Chu, Tosatti, and Weinkove]{CTW}
Jianchun Chu, Valentino Tosatti, and Ben Weinkove.
\newblock On the {$C^{1,1}$} regularity of geodesics in the space of {K\"ahler}
  metrics.
\newblock \emph{Ann. PDE}, 3\penalty0 (2):\penalty0 Paper No. 15, 2017.
\newblock \doi{10.1007/s40818-017-0034-8}.

\bibitem[Crandall et~al.(1992)Crandall, Ishii, and Lions]{CIL}
Michael~G. Crandall, Hitoshi Ishii, and Pierre-Louis Lions.
\newblock User's guide to viscosity solutions of second order partial
  differential equations.
\newblock \emph{Bull. Amer. Math. Soc. (N.S.)}, 27\penalty0 (1):\penalty0
  1--67, 1992.
\newblock \doi{10.1090/S0273-0979-1992-00266-5}.

\bibitem{DR}
T. Darvas and Y. A. Rubinstein,
\emph{Kiselman's principle, the Dirichlet problem for the
Monge--Amp\`ere equation, and rooftop obstacle problems},
J. Math. Soc. Japan \textbf{68} (2016), no.~2,
773--796.

\bibitem[Dinew and Ko{\l}odziej(2017)]{DK}
S{\l}awomir Dinew and S{\l}awomir Ko{\l}odziej.
\newblock Liouville and {Calabi--Yau} type theorems for complex {Hessian}
  equations.
\newblock \emph{Amer. J. Math.}, 139\penalty0 (2):\penalty0 403--415, 2017.
\newblock \doi{10.1353/ajm.2017.0009}.

\bibitem[Donaldson(1999)]{Donaldson}
Simon~K. Donaldson.
\newblock Symmetric spaces, {K\"ahler} geometry and {Hamiltonian} dynamics.
\newblock In \emph{Northern California Symplectic Geometry Seminar}, volume 196
  of \emph{Amer. Math. Soc. Transl. Ser. 2}, pages 13--33. American
  Mathematical Society, Providence, RI, 1999.

\bibitem[Donaldson(2010)]{DonaldsonVolume}
Simon~K. Donaldson.
\newblock Nahm's equations and free-boundary problems.
\newblock In \emph{The Many Facets of Geometry}, pages 71--91. Oxford
  University Press, Oxford, 2010.

\bibitem[Evans(1982)]{Evans}
Lawrence~C. Evans.
\newblock Classical solutions of fully nonlinear, convex, second-order elliptic
  equations.
\newblock \emph{Comm. Pure Appl. Math.}, 35\penalty0 (3):\penalty0 333--363,
  1982.
\newblock \doi{10.1002/cpa.3160350303}.

\bibitem[G{\aa}rding(1959)]{Garding}
Lars G{\aa}rding.
\newblock An inequality for hyperbolic polynomials.
\newblock \emph{J. Math. Mech.}, 8:\penalty0 957--965, 1959.

\bibitem[Ghosh et~al.(2026)Ghosh, Wang, and Xu]{GWX}
Kartick Ghosh, Ziyue Wang, and Yulun Xu.
\newblock Geodesics for generalised complex Hessian equations.
\newblock Preprint, \href{https://arxiv.org/abs/2611.25768}{arXiv:2611.25768}, 2026.

\bibitem[Gilbarg and Trudinger(2001)]{GT}
David Gilbarg and Neil~S. Trudinger.
\newblock \emph{Elliptic Partial Differential Equations of Second Order}.
\newblock Classics in Mathematics. Springer-Verlag, Berlin, 2001.
\newblock Reprint of the 1998 edition.

\bibitem{GP}
P. Guan and D. H. Phong,
\emph{Partial Legendre transforms of non-linear equations},
Proc. Amer. Math. Soc. \textbf{140} (2012), no.~11,
3831--3842.

\bibitem[Guo et~al.(2023)Guo, Smith, and Tong]{GST}
Bin Guo, Kevin Smith, and Freid Tong.
\newblock On the {Hessian-cscK} equations.
\newblock \emph{Math. Z.}, 303\penalty0 (2):\penalty0 Paper No. 33, 2023.
\newblock \doi{10.1007/s00209-022-03187-1}.

\bibitem{GS}
M. J. Gursky and J. Streets,
\emph{A formal Riemannian structure on conformal classes
and uniqueness for the $\sigma_2$-Yamabe problem},
Geom. Topol. \textbf{22} (2018), no.~6, 3501--3573.


\bibitem{HXZ}
W. He, L. Xu, and M. Zhang,
\emph{A fully nonlinear partial differential equation
and its application to the $\sigma_k$-Yamabe problem},
J. Funct. Anal. \textbf{281} (2021), no.~7, Article 109140.

\bibitem[Hou et~al.(2010)Hou, Ma, and Wu]{HMW}
Zuoliang Hou, Xinan Ma, and Damin Wu.
\newblock A second order estimate for complex {Hessian} equations on a compact
  {K\"ahler} manifold.
\newblock \emph{Math. Res. Lett.}, 17\penalty0 (3):\penalty0 547--561, 2010.
\newblock \doi{10.4310/MRL.2010.v17.n3.a12}.

\bibitem[Krylov(1982)]{Krylov}
Nikolai~V. Krylov.
\newblock Boundedly inhomogeneous elliptic and parabolic equations.
\newblock \emph{Izv. Akad. Nauk SSSR Ser. Mat.}, 46\penalty0 (3):\penalty0
  487--523, 1982.
\newblock In Russian.

\bibitem[Mabuchi(1987)]{Mabuchi}
Toshiki Mabuchi.
\newblock Some symplectic geometry on compact {K\"ahler} manifolds. {I}.
\newblock \emph{Osaka J. Math.}, 24\penalty0 (2):\penalty0 227--252, 1987.

\bibitem[Semmes(1992)]{Semmes}
Stephen Semmes.
\newblock Complex {Monge--Amp\`ere} and symplectic manifolds.
\newblock \emph{Amer. J. Math.}, 114\penalty0 (3):\penalty0 495--550, 1992.

\bibitem[Silvestre and Sirakov(2014)]{SS}
Luis Silvestre and Boyan Sirakov.
\newblock Boundary regularity for viscosity solutions of fully nonlinear
  elliptic equations.
\newblock \emph{Comm. Partial Differential Equations}, 39\penalty0
  (9):\penalty0 1694--1717, 2014.
\newblock \doi{10.1080/03605302.2013.842249}.





\end{thebibliography}
\end{document}